\documentclass[aos,seceqn,dvips]{arximspdf}
\usepackage{mathbh}
\usepackage{graphicx}

\doi{10.1214/10-AOS837}
\volume{39}
\issue{1}
\pubyear{2011}
\firstpage{232}
\lastpage{260}

\makeatletter
\newtheorem{theorem}{Theorem}[section]
\newtheorem{lemma}{Lemma}[section]

\newcommand\N{\mathbb{N}}
\newcommand{\EE}{\mathbb{E}}
\newcommand{\RR}{\mathbb{R}}
\newcommand{\HH}{\mathbb{H}}
\newcommand{\R}{\mathbb{R}}
\newcommand{\la}{{\lambda}}
\newcommand{\calD}{{\mathcal D}}
\newcommand{\calF}{{\mathcal F}}
\newcommand{\calH}{{\mathcal H}}
\newcommand{\Var}{\operatorname{Var}}

\newcommand{\calA}{{\mathcal A}}
\newcommand{\calP}{{\mathcal P}}
\newcommand{\calS}{{\mathcal S}}
\newcommand{\calV}{{\mathcal V}}
\newcommand{\calX}{{\mathcal X}}
\newcommand{\Y}{{\mathcal Y}}
\newcommand{\eqd}{\stackrel{{\mathcal D}}{=}}
\newcommand{\toP}{\stackrel{{P}}{\to}}
\newcommand{\X}{{\mathcal X}}
\newcommand{\0}{{\mathbf0}}
\newcommand{\1}{\mathbh{1}}

\newcommand{\tod}{\stackrel{{\mathcal D}}{\longrightarrow}}

\makeatother

\begin{document}
\begin{frontmatter}

\title{Nonparametric estimation of surface integrals}
\runtitle{Estimation of surface integrals}

\begin{aug}
\author[A]{\fnms{Ra\'ul} \snm{Jim\'enez}\thanksref{t1}\ead[label=e1]{raul.jimenez@uc3m.es}}
\and
\author[B]{\fnms{J. E.} \snm{Yukich}\thanksref{t2}\corref{}\ead[label=e2]{joseph.yukich@lehigh.edu}}
\thankstext{t1}{Supported in part by MEC Grant ECO-2008-05080 and CAM
Grant 2008-00059-002.}
\thankstext{t2}{Supported in part by NSF Grant DMS-08-05570.}
\runauthor{R. Jim\'enez and J. E. Yukich}
\affiliation{Universidad Carlos III de Madrid and Lehigh University}
\address[A]{Department of Statistics \\
Universidad Carlos III de Madrid \\
C/Madrid, 126\\ 28903 Getafe (Madrid)\\ Spain\\
\printead{e1}} %adresu isvedimo komanda gale!
\address[B]{Department of Mathematics \\
Lehigh University \\
Bethlehem, Pennsylvania 18015\\
USA \\
\printead{e2}}
\end{aug}

% HISTORY:
\received{\smonth{11} \syear{2009}}
\revised{\smonth{5} \syear{2010}}

% ABSTRACT
%
\begin{abstract}
The estimation of surface integrals on the boundary of an unknown
body is a challenge for nonparametric methods in statistics, with
powerful applications to physics and image analysis, among other
fields. Provided that one can determine whether random shots hit the
body, Cuevas et al. [\textit{Ann. Statist.} \textbf{35} (2007) 1031--1051]
estimate the boundary measure
(the boundary length for planar sets and the surface area for
3-dimensional objects) via the consideration of shots at a box
containing the body. The statistics considered by these authors, as
well as those in subsequent papers, are based on the estimation of
Minkowski content and depend on a smoothing parameter which must be
carefully chosen. For the same sampling scheme, we introduce a new
approach which bypasses this issue, providing strongly consistent
estimators of both the boundary measure and the surface
integrals of scalar functions, provided one can collect the function
values at the sample points. Examples arise in experiments in which
the density of the body can be measured by physical properties of
the impacts, or in situations where such quantities as temperature
and humidity are observed by randomly distributed sensors. Our
method is based on random Delaunay triangulations and involves a
simple procedure for surface reconstruction from a dense cloud of
points inside and outside the body. We obtain basic asymptotics of
the estimator, perform simulations and discuss, via Google Earth's
data, an application to the image analysis of the Aral Sea coast and
its cliffs.
\end{abstract}

% KEYWORDS
%
\begin{keyword}[class=AMS]
\kwd[Primary ]{62G05}
\kwd[; secondary ]{60D05}.
\end{keyword}
\begin{keyword}
\kwd{Surface estimation}
\kwd{boundary measure}
\kwd{Delaunay triangulation}
\kwd{stabilization methods}.
\end{keyword}

\end{frontmatter}

%s1 ###
\section{Introduction}

The estimation of functionals defined on the boundary $\Gamma$ of
an unknown body $G\subset\RR^d$ is a new branch of nonparametric
statistics with powerful applications in several areas, including
image analysis and stereology \cite{cuevas2009}. Cuevas et al. \cite
{cuevas2007} address the estimation of the
Minkowski content
%e1 ###
%
\begin{equation}\label{Mink}
\lim_{\varepsilon\to0^+}
\frac{\mu(\bigcup_{x\in\Gamma} B_\varepsilon(x))}{2\varepsilon},
\end{equation}
$\mu:= \mu_d$ being the Lebesgue measure on $\RR^d$ and
$B_\varepsilon(x)$ the closed $d$-dimensional ball with center at
$x$ and radius $\varepsilon$. When the limit (\ref{Mink}) exists,
it is the most basic measure of the content of $\Gamma$ and it
coincides with length and surface area in dimensions $2$ and $3$,
respectively \cite{mattila1999}. Minkowski content estimators are
based on random point samples distributed on a $d$-dimensional
rectangular solid containing $G$, for which one may determine
whether a point is in $G$ or not. Roughly speaking, they are
empirical measures of the $\varepsilon$-approximation
%e2 ###
%
\begin{equation}\label{Mink2}
\frac{\mu(\bigcup_{x\in\Gamma}
B_\varepsilon(x))}{2\varepsilon}
\end{equation}
of the Minkowski content of $\Gamma$. Both the statistic considered
by Cuevas et al. \cite{cuevas2007} and other closely related statistics
\cite{armendariz2009,cuevas2010,pateiro2008} depend on the
\textit{smoothing parameter}~$\varepsilon$, which must be chosen as a
function of the size of the random point sample.

We propose a different nonparametric approach, free of smoothing
parameters, to estimate not only boundary lengths of planar sets
and surface areas of solids, but also surface integrals
of those scalar functions whose values are knowable at the sample
points. While this paper focuses on instances where the body $G$
is unknown, the proposed method is also of interest for bodies
having a known but complex boundary, one having a surface integral
defying traditional estimators.

The nonparametric estimation of surface integrals has practical
applications in image analysis, when the body $G$ consists of some
nonhomogeneous material and the density $h(x)$ of the material at
$x$ is collected at the sample points. Thus, one might, for example,
be interested in the mass per unit thickness of $\Gamma$, which
corresponds to the surface integral commonly represented by
$\int_\Gamma h\,d\Gamma$. In some instances, quantities such as
temperature and humidity can be measured by sensors randomly
distributed on a set, in which an unknown body is embedded, and a
fundamental problem is to estimate the surface integral of these
quantities on the boundary of the body. In other situations,
including those arising in medical imaging, oncology and cardiology,
knowledge of boundary length is of importance in the prognosis of an
infarction, as well as in the assessment of the dissemination
capacity of a tumor, as explained in Cuevas et al.
\cite{cuevas2007}. %We address this type of problem in this paper.

Our statistics are based on the Delaunay triangulation of the sample points.
Although the idea can be easily adapted to other graphs, in particular
to the Voronoi diagram,
we chose the Delaunay triangulation for two specific reasons:

\begin{enumerate}
\item It is a well-known tool in curve reconstruction methods.
In particular, the $\beta$-skeleton \cite{kirkpatrick1998}, the
Crust \cite{amenta1998} and a wide range of related algorithms
involve computing Delaunay triangulations.
\item
It is computationally efficient---for example, Boissonnat and
Cazals \cite{boissonnat2000} report a 3-dimensional Delaunay
triangulation code which can handle 500,000 randomly distributed
points per minute.
\end{enumerate}

The basic difference between previous methods and the one introduced
here is
that formerly used methods estimate a fattened boundary, namely
$\bigcup_{x\in\Gamma} B_\varepsilon(x)$, whereas we directly estimate
$\Gamma$ by a surface (a curve if $d=2$) properly selected among the
polyhedra (polygons if $d=2$) whose faces belong to the Delaunay
triangulation. Part of our methodology involves a new algorithm for
surface reconstruction based on inner and outer sampling points,
which differs substantially from the numerical methods for surface
reconstruction in which the sample points are on (or close to) the
surface to be reconstructed. Our method is described in Section \ref{sec2},
where we introduce the relevant
statistics.
Section \ref{sec3} establishes basic asymptotics of these statistics and
provides a strongly consistent estimator of the surface integral. In
Section \ref{simulations}, we perform a simulation study. In particular, we estimate
the Minkowski content of sets, comparing our method with existing
ones. An application to image analysis of the Aral Sea coast and its
cliffs from Google Earth's data is discussed in Section \ref{sec5}. Section
\ref{sec6}
summarizes our conclusions. Our proofs, given in Section~\ref{proofs}, rely on
point process methods, including weak convergence of point processes
and stabilization methods, a tool for establishing general limit
theorems for sums of weakly dependent terms in geometric probability \cite{baryshnikov2005,penrose2007EJP,penrose2007Ber,penrose2003}.

%%%%%%%%%%%%%%%%%%%%%%%%%%%%%
%s2 ###
\section{The method: Sewing boundaries of unknown sets}\label{sec2}

Following the basic assumptions of \cite{cuevas2007,cuevas2010},
we will assume that $G$ is a compact subset of an open and bounded
$d$-dimensional rectangular solid $Q$ and that the closure of the
interior of $G$ has positive $\mu$-measure. The boundary of $G$ will
be denoted by $\Gamma$, that is,
\[
\Gamma:= \{x\dvtx \mbox{ for any } \varepsilon> 0, B_\varepsilon(x)
\cap G \neq\varnothing \mbox{ and } B_\varepsilon(x) \cap G^c \neq
\varnothing\},
\]
with $G^c := Q\setminus G$. It is assumed that the \textit{
%e3 ###
$\mu$-boundary of $G$}, defined by
\begin{equation} \label{mub}
\bigl\{x\dvtx\mbox{ for any } \varepsilon> 0, \mu\bigl(B_\varepsilon(x) \cap G\bigr) >0
\mbox{ and } \mu\bigl(B_\varepsilon(x) \cap G^c\bigr)>0 \bigr\},
\end{equation}
coincides with $\Gamma$. This rules out the existence of
``extremities'' to $G$ having null $d$-dimensional Lebesgue measure.
Thus, if we randomly plot enough points in $Q$, there will be points
inside and outside $G$ close enough to any point on the boundary
$\Gamma$. We assume that $\Gamma$ is a $(d-1)$-rectifiable set. That
is, there exists a countable collection of continuously
differentiable maps $g_i\dvtx\RR^{d-1}\to\RR^d$ such that
%e4 ###
%
\begin{equation}
\calH^{d-1}\Biggl(\Gamma\Big\backslash\bigcup_{i = 1}^{\infty}
g_i(\RR^{d-1})\Biggr) = 0,
\end{equation}
$\calH^{d-1}$ being the $(d-1)$-dimensional Hausdorff measure
\cite{mattila1999}. We will assume that $\Gamma$ has finite
Hausdorff measure and that $\Gamma$ has tangent spaces which are
defined almost everywhere.

As in \cite{cuevas2007,cuevas2010}, the sampling model consists of
$n$ i.i.d. random variables $X_1,\dots, X_n$, uniformly distributed
on $Q$, and $n$ i.i.d. Bernoulli random variables
$\delta_1,\dots,\delta_n$ such that
%e5 ###
%
\begin{equation}\label{deldef}
\delta_k =\cases{
1,&\quad if $X_k\in G$, \cr
0,&\quad if $X_k\notin G$.
}
\end{equation}
In other words, although $G$ is unknown, we know whether a sample point $X_k$ is inside $G$ or not. In addition,
given measurable $h\dvtx Q\to\R$, we assume that we know the values
$h(X_1),\dots,h(X_n)$; in general, $h$ is unknown on its domain,
but we assume that we are able to collect its value at each of the
$n$ sample points. Our goal is to estimate the surface integral of
$h$, formally defined by
%e6 ###
%
\begin{equation}\label{intsurf1}
\int_\Gamma h\,d\Gamma= \int_\Gamma h(\gamma) \calH^{d-1}(d\gamma).
\end{equation}
For sets $\Gamma$ satisfying a stronger notion of  $(d-1)$-rectifiability,
namely if $\Gamma$ is the Lipschitz image of a subset of $\R^d$, the
integral (\ref{intsurf1}) coincides (up to a constant factor)
with the Minkowski  content of  $\Gamma$ when $h = 1$.  Thus the
estimated target quantities in
\cite{armendariz2009,cuevas2007,cuevas2010,pateiro2008} can
be written as a surface integral (\ref{intsurf1}) with $h(\gamma) \equiv
1$. This brings  the issues and problems of
\cite{armendariz2009,cuevas2007,cuevas2010,pateiro2008}
within the compass of this paper.

\begin{figure}

\includegraphics{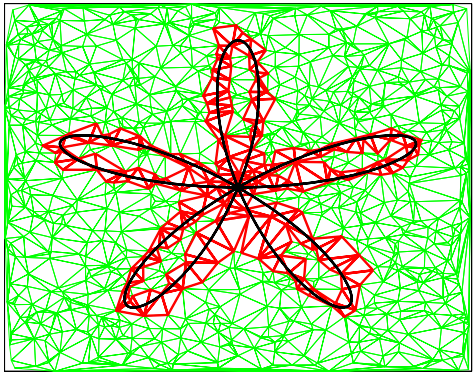}

\caption{Polar rose with polar coordinate equation $\rho= {\frac
{4}{5}}\sin(5\theta)$ (black); Delaunay tessellation for $10^3$
i.i.d. points, uniformly distributed on the square $[-1,1]^2$ (green
and red). Red triangles are the sewing of the rose.}\label{rose_1000}
\end{figure}

Denote by $\calX_n$ the set of sample points $\{X_1,\dots, X_n\}$
and let $\calD(\calX_n)$ be the Delaunay triangulation for
$\calX_n$, namely the full collection of simplices (triangles when
$d=2$) with vertices in $\calX_n$ satisfying the empty sphere
criterion, namely that no point in $\calX_n$ is inside the
circumsphere (circumcircle if $d=2$) of any simplex in
$\calD(\calX_n)$. For any sample of absolutely continuous i.i.d.
points on $\RR^d,$ there exists a unique Delaunay triangulation
almost surely \cite{moller1994}. Each simplex $s \in
\calD(\calX_n)$ is represented by a subset of $(d + 1)$ vertices
belonging to $\calX_n$, here denoted by
$\{X_{s(1)},\dots,X_{s(d+1)}\}$. Recalling (\ref{deldef}), we
introduce the \textit{sewing} of $\Gamma$, denoted by
$S(\calX_n,\Gamma)$, and defined by
\[
\calS(\calX_n,\Gamma):= \Biggl\{s\in \calD(\calX_n)\dvtx 1\leq
\sum_{k=1}^{d+1} \delta_{s(k)} \leq d \Biggr\}.
\]
Thus, $S(\calX_n,\Gamma)$ is the collection of simplices in
$\calD(\calX_n)$ with at least one vertex in $G$ and at least one
vertex in $G^c$. We may assume without loss of generality that there
are no sample points on $\Gamma$. In this case, the sewing of
$\Gamma$ consists of the simplices in the triangulation which
intersect the boundary of $G$. As an illustration, Figure \ref{rose_1000} shows
the sewing of a polar rose with five petals [with polar coordinate
equation $\rho= {\frac{4}{5}}\sin(5\theta)$] based on a sample of
$10^3$ points.
%f1 ###
%

\begin{figure}[b]

\includegraphics{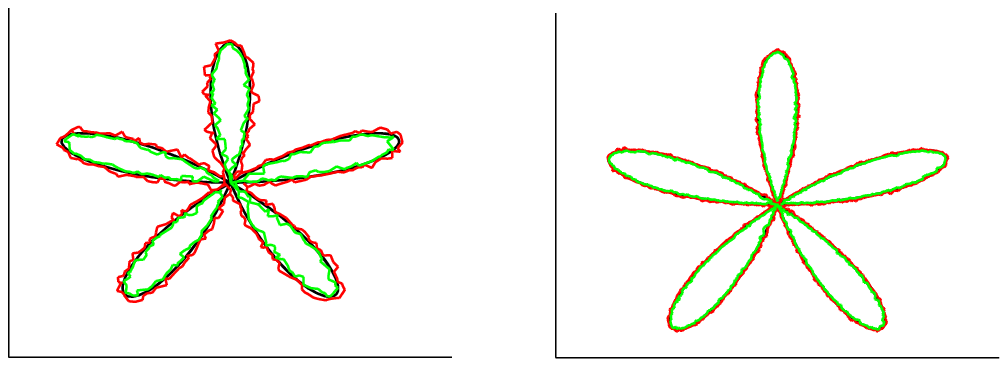}

\caption{Inner (green) and outer (red) sewing of the polar rose (black)
with equation $\rho= {\frac{4}{5}}\sin(5\theta)$, based on
$10^4$ (left) and $10^5$ (right) i.i.d. points uniformly distributed on
the square $[-1,1]^2$.}\label{fig2}
\end{figure}

We are particularly interested in the following two surfaces (curves if $d=2$)
contained in $\calS(\calX_n,\Gamma)$:
\begin{description}
\item[The inner sewing\normalfont{,}] here denoted by $\calS^-(\calX_n,\Gamma),$ which
consists of the union of all faces of the simplices of
$\calS(\calX_n,\Gamma)$ having vertices in $G$.
\item[The outer sewing\normalfont{,}] here denoted by $\calS^+(\calX_n,\Gamma)$,
which consists of the union of all faces of the simplices of
$\calS(\calX_n,\Gamma)$ having vertices in $G^c$.
\end{description}
To be more precise, let $s$ be a simplex of $\calS(\calX_n,\Gamma)$
and $f$ a face of $s$, \textit{where, here and henceforth, by ``face'' we
mean a simplex of dimension $(d-1)$}. Every such face $f$ may be
represented by a vertex set of size $d$, denoted by
%e7 ###
%
\begin{equation}\label{vertexes}
 \calV(f) := \bigl\{X_{f(1)},\dots,X_{f(d)}\bigr\}.
\end{equation}
A face $f$ of $s$ is in the inner sewing $\calS^-(\calX_n,\Gamma)$
if and only if $s\in\calS(\calX_n,\Gamma)$ and $\calV(f) \subset
G$. The face itself need not lie wholly in $G$. On the other hand,
a face $f$ of $s$ is in the outer sewing $\calS^+(\calX_n,\Gamma)$
if and only if $s\in\calS(\calX_n,\Gamma)$ and $\calV(f) \subset
G^c$; again, the face need not lie wholly in $G^c$. Both the inner
and outer sewing consist of polyhedral surfaces (polygons if $d=2$)
which can be used to estimate $\Gamma$. In Figure \ref{fig2}, we show the
inner and outer sewing of the polar rose in Figure \ref{rose_1000}, this
time based on samples with $10^4$ and $10^5$ points. As one would
expect, both
sewings
fit %converge to
the rose when $n$
is large. %\to\infty$.
%f2 ###
%

Now that we have a way to estimate the surface $\Gamma$ by
sampling surfaces, we are ready to explore estimators for the
surface integral at (\ref{intsurf1}). Any numerical approximation
of the integral involving either the inner or the outer sewing is
a potential estimator of the integral on $\Gamma$. In this work,
we approximate integrals by the trapezoidal rule. Thus, if $h$ is
continuous and bounded, the surface integral $ \int_f h\,df$ of any
face $f\in\calS(\calX_n,\Gamma)$ can be approximated by
\[
\int_f h\,df \approx\calH^{d-1}(f)\frac{ 1}{d} \biggl(\sum_{X_k\in\calV
(f)} h(X_k)\biggr) ,
\]
$\calV(f)$ being the vertex set of $f$ defined at (\ref{vertexes}).
Here, $\calH^{d-1}(f)$ is the length of an edge if $d=2$,
the area of a triangle if $d=3$ and the volume of a
$(d-1)$-dimensional simplex otherwise. Therefore, the surface
integrals of $h$ with respect to the inner and outer sewings can be
approximated by the sums
%e8 ###
%
\begin{equation}\label{innerint}
I_{n}^-(h,\Gamma) := \sum_{f\in
\calS^-(\calX_n,\Gamma)} \calH^{d-1}(f) \frac{ 1}{d}
\biggl(\sum_{X_k\in\calV(f)} h(X_k)\biggr)
\end{equation}
and
%e9 ###
%
\begin{equation}\label{outerint}
I_{n}^+(h,\Gamma) := \sum_{f\in
\calS^+(\calX_n,\Gamma)} \calH^{d-1}(f) \frac{ 1}{d}
\biggl(\sum_{X_k\in\calV(f)} h(X_k)\biggr),
\end{equation}
respectively. Next, we study the basic properties of the statistics
(\ref{innerint}) and (\ref{outerint}),
and provide
strongly consistent estimators of $\int_\Gamma h\,d\Gamma$.

%%%%%%%%%%%
%s3 ###
\section{Asymptotics}\label{sec3}

Let $\X\subset\R^d$ be a locally finite point set, that is, for any
compact set $K\subset\R^d$,
$\X\cap K$ contains at most a finite number of points from $\X$.
Let $B \subset\R^d$ be a body with boundary $\partial B$. Define ${\calS}^-(\X, \partial B)$ and   ${\mathcal S}^+(\X, \partial
B)$ analogously to ${\mathcal S}^-(\X_n, \Gamma)$ and
     ${\mathcal S}^+(\X_n, \Gamma)$, respectively. For $x
\in\X$, when $x \in\calS^-(\X, \partial B),$ we let $\xi^-(x, \X,
B)$ be the normalized sum of the Hausdorff measures of the faces
belonging to $\calS^-(\X, \partial B)$ and containing $x$. If $x
\notin\calS^-(\X, \partial B),$ then we put $\xi^-(x, \X, B)= 0.$
Thus,
\[
\xi^-(x, \X, B) :=\cases{
\displaystyle\frac{1}{d}\sum_{\{ f\in{\mathcal S}^-(\X,\partial B)\dvtx x \in{\mathcal V}(f) \} } \calH^{d-1}(f),
&\quad if $ \displaystyle x\in{\bigcup_{f\in\mathcal S^-(\X,\partial B)}} {\mathcal V}(f)$,\hspace*{-2pt}\cr
0, &\quad otherwise.
}
\]
Similarly, we define $\xi^+(x, \X, B)$ to be the normalized sum of
the Hausdorff measures of the faces belonging to $\calS^+(\X,
\partial B)$ and containing $x$; if no such face exists, then $\xi^+$
is defined to be zero.

We will use the functionals $\xi^-$ and $\xi^+$ for several
purposes; in particular, the statistics (\ref{innerint}) and
(\ref{outerint}) can be expressed as the weighted sums
%e10 ###
%
\begin{equation}\label{stat1}
I_{n}^-(h,\Gamma) = \sum_{k=1}^n h(X_k) \xi^-(X_k,\calX_n, G)
\end{equation}
and
%e11 ###
%
\begin{equation}\label{stat2}
I_{n}^+(h,\Gamma) = \sum_{k=1}^n h(X_k) \xi^+(X_k,\calX_n, G).
\end{equation}

The functional $\xi^-$ is translation invariant, that is, for all $x
\in\R^d$, all locally finite $\X$ and all bodies $B$, we have
$\xi^-(x,\X, B) = \xi^-(\0,\X-x, B-x)$; here, $\0$ is a point at the
origin of $\RR^d$ and for sets $F \subset\R^d$ and $x \in\R^d,$ we
put\break $\{F - x\}:=\{y-x\dvtx y \in F\}$. Also, given $\alpha> 0$ and
putting $\alpha F := \{\alpha y\dvtx y \in F\}$, we have $
\calH^{d-1}(\alpha f) = \alpha^{d-1} \calH^{d-1}(f)$. Thus, $\xi^-$
satisfies the following {\em scaling property} for all $\eta>0$:
%e12 ###
%
\begin{equation}
\eta^{d-1}\xi^-(\0,\X, B) = \xi^-(\0,\eta\X,\eta B),
\end{equation}
%
%e13 ###
which, when combined with translation invariance, gives
\begin{equation}\label{scaling}
 \eta^{d-1}\xi^-(x,\X, B) =\xi^-\bigl(\0,\eta(\X-x),\eta(B-x)\bigr).
\end{equation}

Thus, by the definition of $I_{n}^-(h,\Gamma),$ we have
%e14 ###
%
\begin{equation} \label{scaling1}
\hspace*{15pt}I_{n}^-(h,\Gamma) = n^{-(d-1)/d} \sum_{k=1}^n h(X_k)
\xi^-\bigl(\0,n^{1/d}(\calX_n-X_k), n^{1/d}(G-X_k)\bigr).
\end{equation}
Similarly, $\xi^+$ is translation invariant and satisfies the
scaling property (\ref{scaling}) and, therefore,
%e15 ###
%
\begin{equation} \label{scaling2}
\hspace*{15pt}I_{n}^+(h,\Gamma) = n^{-(d-1)/d} \sum_{k=1}^n h(X_k)
\xi^+\bigl(\0,n^{1/d}(\calX_n-X_k), n^{1/d}(G-X_k)\bigr).
\end{equation}

Central to our results is the asymptotic behavior of the summands in
(\ref{scaling1}) and~(\ref{scaling2}); see Lemma \ref{lemma1}. For
this, it is convenient to introduce the random
variable $\xi(t),$ defined at (\ref{defxi}) below. %; see Lemma

Denote by $\calP$ the homogeneous Poisson point process of intensity
1 on $\RR^d$ and let $\calP^{\0} := \calP\cup\{\0\}$. In
general, for any $\lambda>0$, we denote by $\calP_\lambda$ the
homogeneous Poisson process of intensity $\lambda$ on $Q$. %[$\calP_

For all $t \in\R$, denote by $\HH_t^d$ the \textit{half-space}
%e16 ###
%
\begin{equation} \label{half}
\HH_t^d := \RR^{d-1} \times(-\infty, t].
\end{equation}
The homogeneity of the Poisson point process implies, for all $t \in
\R,$ that $\xi^-( \0, \calP^{\0}, \HH^d_t)$ and $\xi^+( \0,
\calP^{\0}, \HH^d_{-t})$ are equally distributed.

Given $t \in\R$, denote by $\xi(t)$ the normalized sum of the
Hausdorff measures of the faces incident to $\0$ belonging to the
inner sewing ${\mathcal S}^-(\calP^{\0},\partial\HH_t^d)$. If there are no such
faces, then $\xi$ is set to be zero, that is,
%e17 ###
%
\begin{equation} \label{defxi}
\xi(t):= \xi^-( \0, \calP^{\0}, \HH^d_t) \eqd
\xi^+( \0, \calP^{\0}, \HH^d_{-t}).
\end{equation}
Lemmas \ref{lemma2} and \ref{lemma4} imply that $\xi(\cdot)$ is
dominated by an integrable function, therefore, the constant
%e18 ###
%
\begin{equation}\label{alphad}
\alpha_d := \int_0^{\infty} \EE[\xi(t)]\,dt,
\end{equation}
which plays an important role in the asymptotics of our method, is
well defined.

It is also meaningful to consider the Poissonized versions of
$I_{n}^-(h,\Gamma)$ and $I_{n}^+(h,\Gamma)$,
namely,
%e19 ###
%
\begin{equation}\label{Pois}
{\mathcal I}_{\la}^-(h,\Gamma) := \sum_{x \in\calP_\lambda} h(x) \xi^-(x,\calP_\lambda, G)
\end{equation}
and
%e20 ###
%
\begin{equation}\label{Pois2}
{\mathcal I}_{\la}^+(h,\Gamma) := \sum_{x \in
\calP_\lambda} h(x) \xi^+(x,\calP_\lambda, G).
\end{equation}

The following theorems are our main technical results. They take
into account the possibility that $h$ can be discontinuous on
$\Gamma$ from either outside or inside $G$, which is not uncommon
in applications, including, for example, the case when $h$ denotes
density of the body $G$. We say that $h$ is \textit{inner continuous} if
the restriction of $h$ to $G$ is continuous; likewise, we say that $h$
is \textit{outer continuous} if its restriction to the closure of $G^c$ is
continuous.

\begin{theorem}\label{expectation}
Let $\Gamma$ be the boundary of a compact set
$G\subset\RR^d$. Assume that $\Gamma$ is a $(d-1)$-rectifiable set,
that it has finite Hausdorff measure
and coincides with the $\mu$-boundary of $G$, defined at
(\ref{mub}). If
$h$ is inner continuous, we have
\[
\lim_{\la\to\infty} \EE[{\mathcal I}_{\la}^-(h,\Gamma) ] = \lim
_{n\to\infty} \EE[I_{n}^-(h,\Gamma) ] = \alpha_d
\int_\Gamma h\,d\Gamma
\]
and
\[
\lim_{\la\to\infty} {\mathcal I}_{\la}^-(h,\Gamma) = \lim_{n\to
\infty} I_{n}^-(h,\Gamma) = \alpha_d
\int_\Gamma h\,d\Gamma\qquad \mathrm{a.s.},
\]
with $\alpha_d$ defined at (\ref{alphad}). If
$h$ is outer continuous,
then we have the same asymptotics for ${\mathcal I}_{\la}^+$ and $I_{n}^+$.
\end{theorem}

\begin{theorem}\label{variance}
 Let $\Gamma$ be as in Theorem \ref{expectation}.
If $h$ is inner
continuous, we have
\[
\lim_{\la\to\infty} \lambda^{(d-1)/d} \Var[{\mathcal I}_{\lambda
}^-(h,\Gamma) ] =
V_d
\int_\Gamma h^2\,d\Gamma,
\]
where
%e21 ###
%
\begin{equation}\label{Vd}
V_d:= \int_0^{\infty} \EE[\xi^2(t)]\,dt + \int_0^{\infty}\!\!\!\int_{\RR^{d} } c_t(z)\,dz\,dt,
\end{equation}
with
%e22 ###
%
\begin{equation}\label{cdef}
c_t(z):= \EE[\xi^-( \0, \calP^\0 \cup\{z\},\HH^d_t)\xi^-( z,
\calP^{\0} \cup\{z\}, \HH^d_t)] - (\EE[\xi(t)])^2.
\end{equation}
If $h$ is outer continuous, we have the same asymptotics for
$\Var[{\mathcal I}_{\lambda}^+(h,\Gamma)]$.
\end{theorem}

We expect that modifications of the stabilization arguments in
[\cite*{baryshnikov2005}, \cite*{penrose2007EJP}--\cite*{penrose2003}]
will yield a de-Poissonized version of Theorem
\ref{variance}, that is, variance asymptotics for
$I_{n}^-(h,\Gamma)$. These involve nontrivial arguments and we
postpone this to a later paper. Roughly speaking,
since Poisson input introduces more variability than binomial
input, we expect the variances of
$I_{n}^-(h,\Gamma)$ and $I_{n}^+(h,\Gamma)$
to be no larger than the variances of $ {\mathcal I}_{n}^-(h,\Gamma)$ and
$ {\mathcal I}_{n}^+(h,\Gamma)$, respectively.
Consequently, we have good reasons to believe that, under the
assumptions of Theorem \ref{variance}, both
$\Var[I_{n}^-(h,\Gamma)]$ and $\Var[I_{n}^+(h,\Gamma)]$ are
$O(n^{-(d-1)/d})$. We have not yet obtained analogous
results for the rate of convergence for the bias
$(\EE[I_{n}^{\pm}(h,\Gamma)]- \alpha_d \int_\Gamma
h\,d\Gamma)$. In this direction, we provide, in Section
\ref{simulations}, some experimental results for dimension $d=2$.

In accordance with the assumptions of Theorems \ref{expectation}
and \ref{variance}, we define
\[
S_n(h,\Gamma):=\cases{
I_{n}^-(h,\Gamma),
& \quad if $h$ is inner continuous on $G$,\cr
&\quad but not outer continuous on $G$,\cr
I_{n}^+(h,\Gamma) , & \quad if $h$ is outer continuous on $G$,\cr
& \quad but not inner continuous on $G$, \cr
\frac{1}{2} \bigl(I_{n}^-(h,\Gamma) +I_{n}^+(h,\Gamma)\bigr),
& \quad if $h$ is continuous everywhere.\label{estimator}
}
\]
In the light of our results and remarks, we propose to estimate the
surface integral $\int_\Gamma h\,d\Gamma$ by the strongly
%e23 ###
consistent sewing-based estimator
\begin{equation}\label{estimator1}
 \alpha_d^{-1}S_n(h,\Gamma).
\end{equation}

Since Theorems \ref{expectation} and \ref{variance} only assume the
rectifiability of $\Gamma$, the estimator (\ref{estimator1})
is applicable when the body $G$ has sharp interior or exterior cusps.
In particular, it can be used for estimating the boundary measure of
such sets, an estimation problem in which previous methods \cite{armendariz2009,cuevas2007,pateiro2008} have drawbacks.
Under suitable conditions on the smoothing parameter, \cite
{cuevas2010} also discusses the consistency of the estimators based on
the Minkowski content for a general class of bodies.

%%%%%%%%%%
%s4 ###
\section{Simulations}\label{simulations}

The estimator (\ref{estimator1}) depends on the constant $\alpha_d$
defined at (\ref{alphad}). This constant can be estimated by
Monte Carlo methods, as follows. Consider large $n$ %(stretching the
and a surface $\Gamma^0$
such that $\int_{\Gamma^0} d\Gamma^0$ is known. Let $\1$ denote
the function identically equal to $1$ and recall the definition of our estimator
$S_n$. Then simulate a random sample of $S_{n}(\1,\Gamma^0)$
and estimate $\alpha_d $ by
\[
\hat{\alpha}_d := \frac{\operatorname{mean}(S_{n}(\1,\Gamma^0) )}{\int_{\Gamma^0}\, d\Gamma^0}\label{alpha}
\]
with ${\rm mean}(S_{n}(h,\Gamma^0))$ being the sample mean. Given
$\hat{\alpha}_d$ and an arbitrary surface $\Gamma$, the natural
estimator of $\int_\Gamma h\,d\Gamma$ is thus
%e24 ###
%
\begin{equation}\label{estimator2}
\hat{\alpha}_d^{-1}S_n(h,\Gamma).
\end{equation}

Taking this procedure into account, we estimated $\alpha_2$ by
performing a simulation with $1000$ independent copies of
$S_{n}(\1,\Gamma^0)$, with $\Gamma^0$ being the unit circle and
with sewings based on $n = 10^7$ points
uniformly distributed on the square $[-2,2]^2$.
For this configuration, we obtained
%e25 ###
%
\begin{equation} \label{defalpha}
\hat{\alpha}_2 = 1.1820.
\end{equation}

To compare our estimator with those based on empirical
approximation of the Minkowski content, we estimated the length $L
=12\sqrt{3}\approx20.7846$ of ${\mathcal T}$, ${\mathcal T}$ being
Catalan's trisectrix (also called the Tschirnhausen cubic), with polar
equations
\[
r=\cases{
\sec^3(\theta/3), &\quad if $0 < \theta\leq\pi$, \cr
\sec^3\bigl((2\pi-\theta)/3\bigr),&\quad if $\pi< \theta\leq2\pi$.
}
\]
We used the same simulation framework as used in \cite{cuevas2007},
where $n$ points ($n=10{,}000$ and $n=30{,}000$)
were drawn from the square $[-9,2]\times[-5.5, 5.5]$ $500$ times. As
mentioned above, the estimator proposed in \cite{cuevas2007}
depends on the smoothing parameter $\varepsilon,$ which must be
carefully chosen. Of the fourteen possible values of the smoothing
parameter considered in \cite{cuevas2007} and of the fourteen
different corresponding estimators, we let $L_n^{M}({\mathcal T})$ be
the Minkowski content estimator having, on average, the smallest bias,
that is, the one minimizing the difference between the expectation of
the estimator
and its target value. In Table \ref{table1}, we compare $L_n^{M}({\mathcal T})$ with
the \textit{sewing-based estimator} $L_n^S ({\mathcal T})$, namely
%e26 ###
%
\begin{equation}\label{sewing-based}
 L_n^S ({\mathcal T}) =
\hat{\alpha}_d^{-1}S_n(\1,{\mathcal T}).
\end{equation}
%
%t1 ###
%
\begin{table}\tablewidth=310pt
\caption{Mean and standard deviation of estimators based on the
Minkowski content, with optimal smoothing parameter, and the
sewing-based estimators
over 500 replications and simple sizes $n = 10{,}000$ and $n=30{,}000$
(the true value is $12\sqrt{3}\approx20.7846$)}\label{table1}
\begin{tabular*}{310pt}{@{\extracolsep{\fill}}lcccc@{}}
\hline
$\bolds{n}$ &\textbf{mean$\bolds{(L_n^{M}({\mathcal T}))}$} &$\bolds{\operatorname{std}(L_n^{M}({\mathcal
T}))}$
&\textbf{mean$\bolds{( L_n^{S}({\mathcal T}))}$}&$\bolds{\operatorname{std}(L_n^{S}({\mathcal T}))} $\\
\hline
$10{,}000$&19.3679&1.0394&20.7030& 0.3140\\
$30{,}000$&19.8237&1.0666&20.7328& 0.2375\\
\hline
\end{tabular*}
\end{table}
Roughly speaking, we improve the mean relative error (difference
between exact and approximate mean value/exact value) by a factor
of
$10^{-3}$ and we significantly reduce the standard
deviation.
\begin{figure}[b]

\includegraphics{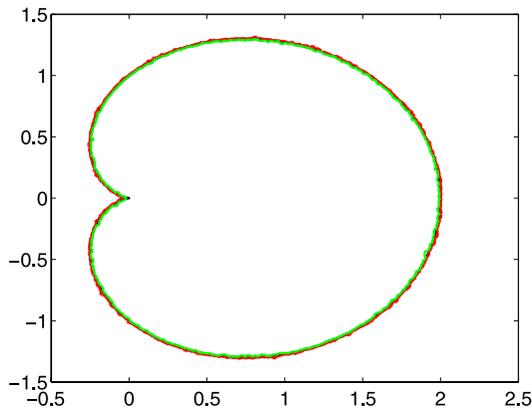}

\caption{Inner (green) and outer (outer) sewing based on $10^5$
uniform points
overlapping a cardioid.}\label{cardioid}
\end{figure}

\begin{figure}[b]

\includegraphics{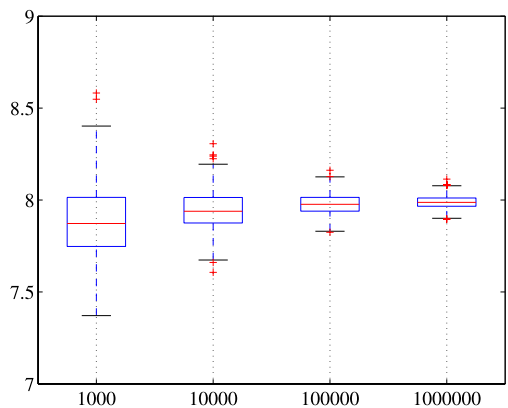}

\caption{Box plots of $1000$ replications of the sewing-based estimator
$L_n^S({\mathcal C})$,
for $n=10^3$, $10^4$, $10^5$ and $10^6$.}\label{boxplot}
\end{figure}

Now, consider the planar cardioid $\mathcal C$ with polar equation $\rho
= 1 +\cos\theta$, $0 \leq\theta\leq2 \pi.$ We selected this curve
since the sharp cusp at the origin precludes using methods based on
empirical approximation of the Minkowski content. Figure \ref{cardioid} shows how
the inner and outer sewing practically overlap the curve with
samples of size $10^5$ of randomly distributed points on the square
$[-0.5,2.5]\times[-1.5,1.5]$.
The target now is to estimate the length of the cardioid, which
equals 8, using the sewing-based estimator $L_n^S({\mathcal C})$. We
performed $1000$ independent copies of $L_n^S({\mathcal C})$ for
$n=10^3$, $10^4$, $10^5$ and $10^6$. The results are summarized in Table
\ref{table2} and Figure~\ref{boxplot}.
%t2 ###
%
\begin{table}\tablewidth=285pt
\caption{Mean, standard deviation and $\sqrt{n}$-times mean square error
of the sewing-based estimator $L_n^S({\mathcal C})$, computed over $1000$
independent replications}\label{table2}
\begin{tabular*}{285pt}{@{\extracolsep{\fill}}lccc@{}}
\hline
$\bolds{n}$ &\textbf{mean}$\bolds{(L_n^S({\mathcal C}))}$ & $\bolds{\operatorname{std}(L_n^S({\mathcal C}))}$&
$\bolds{\sqrt{n}\times}$\textbf{MSE}$\bolds{(L_n^S({\mathcal C}))}$\\
\hline
\phantom{000,,,}$1000$&7.8862&0.1870&1.5157\\
\phantom{00$,$}$10{,}000$&7.9446&0.1009&1.3254\\
\phantom{0$,$}$100{,}000$&7.9772&0.0538&1.0801\\
$1{,}000{,}000$&7.9885&0.0323&1.1736\\
\hline
\end{tabular*}
\end{table}
%
%f4 ###
%
The following important remarks can be drawn from this case study:
\begin{enumerate}
\item As one might expect, the sharp inner cusp is slightly
underestimated. The gain achieved by
increasing the sample size is for both bias and
standard deviation.
\item The simulations strongly suggest that $L_n^S({\mathcal C})$ is very
close to normal.
This result can be used to provide confidence intervals for the
integral that we are estimating. For this, in real applications
where there are no replicas of the estimator, bootstrap procedures
can be used to estimate $\Var[S_n(h,\Gamma)]$.
\item %As we have discussed in previous section,
The mean square error scales with $\sqrt{n}$. %$n^{-(d-1)/d}$.
Thus, the sewing-based estimator seems to converge much faster than
estimators based on the empirical approximation of the Minkowski
content, which, in general, we do not expect to converge faster than
$n^{1/4}$ \cite{armendariz2009}.

\end{enumerate}

We finalize our simulation study by estimating line integrals of a
parametric family $\{h_\zeta, \zeta\in\RR\}$ of scalar functions on
the cardioid ${\mathcal C}$ given by $\rho= 1 + \cos\theta, 0 \leq
\theta\leq2 \pi$.
We choose $h_\zeta(x,y) =\rho^{\zeta+ 1/2}$
%with $\rho= \sqrt{x^2+y^2}$,
for two reasons:
\begin{enumerate}
\item[(i)]
we have simple expressions for the integrals
\[
I(\zeta) = \int_{\mathcal C} h_\zeta\,d{\mathcal C} = \sqrt{2} \int
_0^{2\pi}(1+\cos\theta)^{1 + \zeta}\,d\theta,
\]
where we use $d{\mathcal C} = \sqrt{ (\rho(\theta))^2 +
(\rho'(\theta))^2 }\,d \theta= \sqrt{2\rho(\theta)}\,d \theta$;

\item[(ii)]
the contribution to the total integral of the integral near the
sharp inner cusp increases with $\zeta$, and the variability of
$h_\zeta$ also increases with $\zeta$.
\end{enumerate}

As before, we considered $1000$ independent replications of the
sewing-based estimator of $I(\zeta)$ based on $n=10^3$, $10^4$, $10^5$,
$10^6$ points uniformly distributed on the square $[-0.5,2.5]\times
[-1.5,1.5]$. Similarly to the previous case study, increased sample
size resulted in smaller bias and smaller standard deviation. %the
Also, the simulations
strongly suggested that the estimator is very close to normal and that
the mean square error scales with $\sqrt{n}$.
Our
goal now is to highlight how the estimation depends on the parameter
$\zeta$. For this, we summarize the results for $n=10^4$ in Table \ref{table3}.
%t3 ###
%
\begin{table}\tablewidth=285pt
\caption{Mean and standard deviation, computed over $1000$ independent
replications,
of the sewing-based estimator of $I(\zeta),$ based on
$n=10^4$ points uniformly distributed on the square
$[-0.5,2.5]\times[-1.5,1.5]$}\label{table3}
\begin{tabular*}{285pt}{@{\extracolsep{\fill}}lccc@{}}
\hline
$\bolds{\zeta}$ &$\bolds{I(\zeta)}$&\textbf{mean}$\bolds{(\hat{I}(\zeta))}$ & $\bolds{\operatorname{std}(\hat{I}(\zeta))}$\\
\hline
$-1/2$&\phantom{0}8.0000&\phantom{0}7.9446&0.1009\\
\phantom{$-$}$0$&\phantom{0}8.8858&\phantom{0}8.8786&0.1138\\
\phantom{$-$}$1$&13.3286&13.3349&0.1827\\
\phantom{$-$}$2$&22.2144&22.2191&0.3470\\
\hline
\end{tabular*}
\end{table}
This behavior with respect to $\zeta$ is common for all of the $n$ values that we
considered. We remark that, on the one hand, the bias
is slightly smaller, to the extent that the sharp inner cusp %affects
%weighs
has less of an effect on the value of the integral. On the other hand,
the more
the function $h$ differs from being a constant (i.e., the case
$\zeta= -1/2$), the more the variance of $\hat{I}(\zeta)$ grows.

%
%s5 ###
\section{Application to image analysis using Google Earth data}\label{sec5}

Our simulations show that we require about $10^5$ points or more to
get attractive estimators of the full image of a planar set. But
they also show that 1000 points are enough to obtain good estimates
of integrals along curves. This sample size is quite manageable for
online applications for any new-generation personal computer, on
which a compiled version of our algorithm would run in less than a
second. Even $n=10^4$ could be implemented for snapshot queries.
Next, we discuss an application to image analysis using Google
Earth data. Quite possibly, Google provides a sharp image of the
area around your residence or workplace, but what can you infer
from the image of the most inaccessible regions of the planet,
including, for example, the coast of the Aral Sea?

According to the United Nations, the disappearance of the Aral Sea,
once the fourth largest inland body of water on the planet, is the
worst man-made environmental disaster of the 20th century. It is
estimated that the size of the Aral Sea has fallen by more than 60\%
since the 1960s. Water withdrawal for irrigation has caused a
dramatic fall in the water level, revealing a fascinating
geology of cliffs overlooking the Aral Sea. %Internet shows
What is the mean height of the cliffs along a long waterfront? How
irregular are they?
%Google Earth only displays the big picture of the Aral Sea, any zoom
%in it shows a fuzzy image without traces of cliffs. But,
Google Earth provides elevation information for any pixel and, from
its color, we can determine whether it is in the Aral Sea or not.
This scenario fits our sampling model and thus we may use sewings
to estimate the mean and variance of the elevations of cliffs along
waterfronts. Both parameters are related to common topographic
measures. According to \cite{evans1972}, the standard deviation of
elevation provides one of the most stable measures of the vertical
variability of a topographic surface. On the other hand, the
elevation mean is related to the elevation-relief ratio, namely
\[
\frac{\mbox{elevation mean -- elevation min.}}{\mbox{elevation max. --
elevation min.}},
\]
which has been computed in the past using a point-sampling
technique, rather than planimetry \cite{pike1971}.

\begin{figure}

\includegraphics{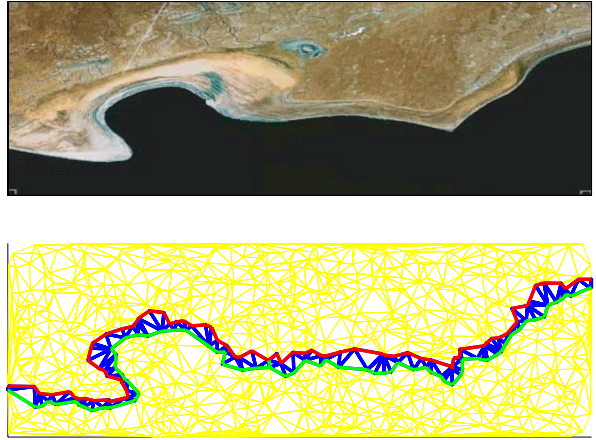}

\caption{Google Earth image of the land area bounded between
latitudes 45.9 and 46.04, and longitudes 58.9 and 59.27 (top);
Delaunay triangulation for 1000 points uniformly distributed on the
rectangle $[58.9,59.27]\times[45.9, 46.04]$ and sewings of the Aral
Sea (bottom).}\label{aral}
\end{figure}

We focus our analysis on the land area bounded between latitudes
45.9 and 46.04, and longitudes 58.9 and 59.27. We chose this surface,
approximately 445~km$^2$, partly because of its complex shape, involving
capes and bays.
%that can be seen as
%smooth versions of inner and outer cusps discussed in previous section.
Figure~\ref{aral}~(top) shows the Google Earth image of the surface under consideration.

We emphasize the three following important points:
\begin{enumerate}
\item
This surface is not a planar rectangle, but we proceed as though
it were. The land area that concerns us is relatively small and
therefore our statistics do not depend on whether we use
Euclidean or geodesic distance. The error produced by considering
land areas as planar sets is briefly discussed later.
\item
One of our basic assumptions is that the set to be estimated is
compact and contained in the interior of the rectangle to be
scanned. This is not the case now because the boundary between sea
and land is not a loop contained by the scanned rectangle. In these
cases, boundary effects can induce
spurious long faces in the extremities of the inner and outer
sewings. It is appropriate to exclude these two faces from each
sewing. We remark that this procedure does not affect the result for
large $n$.
\item The fractal-like structure of coast lines is not detectable with
the limited resolution of Google Earth and thus we may safely assume
that the Aral Sea coast, even if naturally fractal-like, is a
rectifiable curve.
\end{enumerate}
%f5 ###
%

To summarize and clarify, we use the following procedure: drop $1000$~points,
uniformly distributed, on the longitude/latitude rectangle $Q
= [58.9,59.27]\times[45.9, 46.04]$, obtain their corresponding
elevation from Google Earth, determine whether they are in the sea
or not and, finally, compute the sewings of the sample. The results
are shown in Figure \ref{aral} (bottom). In this context, distances between
points are found as follows.
Setting the Earth's radius to be $6371$ km, we approximate
the distance in kilometers between two latitude/longitude points
$(\mbox{long}_1, \mbox{lat}_1)$ and $(\mbox{long}_2, \mbox{lat}_2)$
of $Q$ by the spherical law of cosines:
\[
6371 [\cos^{-1}(\sin(\mbox{lat}_1) \sin(\mbox{lat}_2)) +
\cos(\mbox{lat}_1) \cos(\mbox{lat}_2) \cos(\mbox{long}_2 - \mbox{
long}_1 )].
\]
Alternatively, we obtain the same numerical results using the
appropriately scaled Euclidean distance or the Haversine formula.

Let $\Gamma$ denote the waterfront, that is, the curve contained in
$Q$ separating land from water in the image, and let $L(\Gamma)$
denote the
length of $\Gamma$. Let $G$ represent the land and $h(x,y)$ denote the
height in
meters from the water of the latitude/longitude point $(x,y)\in G$.
We are here assuming that the discrete data given by Google Earth can be
extended to yield a continuous approximation of height at those
points where there is not enough available information. It is reasonable to assume that
the restriction of $h$  to $G$ is continuous, namely  that $h$ is inner
continuous on $G$.

We estimate the mean of the height of the cliffs along $\Gamma$,
namely
\[
\bar{h} := \int_\Gamma h \frac{d\Gamma}{L(\Gamma)},
\]
and the standard deviation
\[
s_{h} := \sqrt{\int_\Gamma(h - \bar{h})^2
\frac{d\Gamma}{L(\Gamma)}},
\]
by
\[
\frac{I^-(h,\Gamma)}{I^-(\1,\Gamma)}= 10.2787\mbox{ m}\quad \mbox{and}
\quad\sqrt{\frac{I^-( (h-\bar{h})^2,\Gamma)}{I^-(\1,\Gamma)}} =
9.8102\mbox{ m},
\]
respectively. Note that $\bar{h}\approx s_{h} $, which certainly
means a significant variance of the height of cliffs along the
waterfront.

We remark that one can compute the above estimators without knowing
the value of $\alpha_d$. Using our estimation of this constant, the
length of the coastline $L(\Gamma)$ is estimated using
(\ref{defalpha}) and (\ref{sewing-based}), yielding
\begin{eqnarray*}
L^{{\mathcal S}}(\Gamma) &=& \frac{I^-(\1,\Gamma) +I^+(\1,\Gamma)}
{2\hat{\alpha}_2} \\
&=& 40.977\mbox{ km}.
\end{eqnarray*}
This length may be contrasted with the distance between the lower
vertices of $Q$, which is approximately $28.6313$ km. All of the above
estimates could be easily implemented using Google Earth.
This requires that the user enter the southwest and northeast
coordinates, and
that there is a rule (the pixel color, in our case)
characterizing the set to be estimated. Of course, the larger the land
area, the greater the differences between spherical and Euclidean
measures. Thus, the larger the area, the worse the estimate. We should
mention that more sophisticated applications, with which we could
provide approximations of the estimation errors, using re-sampling,
for example, could still be costly for online consultation with
current processors.

%
%s6 ###
\section{Conclusions}\label{sec6}

Random scanning of unknown bodies is a good alternative to regular scanning
when the underlying morphology is complex.
The paper of Cuevas et al. \cite{cuevas2007} addressed the problem of
estimating the boundary measure of a set for the former type of these
scanning setups.
For this sampling scheme, we introduce an
efficient computational method for estimating not only
the boundary measure, but also surface integrals of scalar
functions, provided one can collect the function values at the sample
points.

We discuss conditions for getting strong consistency, as well as some
issues related to the rate of convergence of our estimators. Our
proofs rely on point process methods, including weak convergence of
point processes and stabilization methods.

We perform a simulation study to compare our estimators with
previous estimators of boundary lengths of sets, %that satisfy a kind
%of smoothness condition, namely the `free rolling condition'. We
concluding that the sewing-based estimators (\ref{estimator2})
significantly reduce the errors and computation times, while
increasing precision. We complete our simulation study by
estimating boundary lengths of sets with sharp cusps, as well as the
integrals of scalar functions on the boundaries of these sets, always
obtaining good results.

An online application to image analysis using Google Earth data is
discussed. In particular, a complex waterfront of the Aral Sea,
approximately 41 km according to our estimators, is analyzed.
Specifically, we estimate surface integrals related to the mean and
standard deviation of the height of the irregular cliffs facing the
sea line.

%s7 ###
\section{Proofs and technical results}\label{proofs}

As indicated at the outset, our approach makes use of stabilizing
functionals $\xi$, where $\xi(x, \X)$ is a translation invariant
functional defined on pairs $(x, \X)$, where $x \in\RR^d$
and $\X\subset\RR^d$ is locally finite. Translation invariance
means that $\xi(x, \X) = \xi(x +y, x + \X)$ for all $y \in\RR^d$.
We recall a few facts about such functionals [\cite*{baryshnikov2005},
 \cite*{penrose2007EJP}--\cite*{penrose2003}]. As in
\cite{penrose2007EJP}, we consider the following metric on the space
$\mathcal L$ of locally finite subsets of $\R^d$:
%e27 ###
%
\begin{equation} \label{metric}
D(\calA, \calA') := \bigl(\max\{K \in\N\dvtx \calA\cap B_K(\0) = \calA'
\cap B_K(\0)\}\bigr)^{-1}.
\end{equation}

We say that $\xi(\cdot, \cdot)$ stabilizes on a homogeneous Poisson
point process $\calP$ if, for all $x \in\R^d$, there is an a.s.
finite $R:=R(x,\calP)$ such that
\[
\xi\bigl(x, \calP\cap B_R(x)\bigr) = \xi\bigl(x, [\calP\cap B_R(x)] \cup{\mathcal A}\bigr)
\]
for all ${\mathcal A} \subset(B_R(x))^c$. Recall that  $\mathcal P^{\0}:= \mathcal P \cup \{\0\}$. Now, whenever $\xi(\cdot,
\cdot)$ satisfies stabilization, then $\calP^\0$ is a continuity
point for the function $g(\calA):= \xi(\0,\calA)$ with respect to
the topology on $\mathcal L$ induced by $D$. Thus, by the continuous
mapping theorem (Proposition A2.3V on page 394 of \cite{daley2003}), if
${\mathcal Y}_n$, $n \geq1,$ is a sequence of random point measures
satisfying ${\mathcal Y}_n \tod \mathcal{P}^\0$ as $n \to\infty$, then
$\xi(\0,{\mathcal Y}_n) \tod\xi(\0,\calP^\0)$; see
\cite{penrose2007EJP} for details.

Thus, if $U_k$, $k \geq1,$ are i.i.d.
uniform on the unit cube with ${\mathcal U}_n := \{U_k\}_{k=1}^n$ and
if $\xi$ is translation invariant so that $\xi(n^{1/d}U_1, n^{1/d}
{\mathcal U}_n) = \xi(\0, n^{1/d}({\mathcal U}_n - U_1))$, then, since the
shifted and dilated random point measures $n^{1/d}({\mathcal U}_n -
U_1)$ satisfy the convergence $n^{1/d}({\mathcal U}_n - U_1) \tod \mathcal{P}^\0$,
it follows from stabilization that, as $n \to\infty$,
%e28 ###
%
\begin{equation} \label{JY1}
\xi(n^{1/d}U_k, n^{1/d} {\mathcal U}_n) \tod\xi(\0, \calP^\0).
\vadjust{\goodbreak}\end{equation}

This result is central to proving the weak law of large numbers
\cite{penrose2003},
%e29 ###
%
\begin{equation} \label{weaklaw}
n^{-1}\sum_{k=1}^n \xi(n^{1/d}U_k, n^{1/d}{\mathcal U}_n) \toP
\EE[\xi(\0,\calP^\0)].
\end{equation}
The analogous convergence result for the two-dimensional
vector
\[
[\xi(n^{1/d}U_k, n^{1/d} {\mathcal
U}_n), \xi(n^{1/d}U_j, n^{1/d} {\mathcal U}_n)],\qquad k \neq j,
\]
is likewise
central to showing variance asymptotics and central limit theorems
\cite{baryshnikov2005,penrose2007EJP} for the sums in
(\ref{weaklaw}).

One of the main features of our proofs is that if $\xi$ is
stabilizing, if $\Y_n$, $n \geq1,$ is a sequence of point
processes on sets increasing up to the half-space $H$ which is
a translation/rotation of $\HH_0^d := \RR^{d-1}\times(-\infty,0]$
and if $\Y_n \tod\calP\cap H$ as $n \to\infty$, then, subject to
the proper scaling, there exist limit results for $\xi$ analogous to
those at (\ref{weaklaw}). Our goal is to make these ideas precise.

\textit{We prove our main results for the functional $\xi^-$.
Identical results hold for the
functional $\xi^+$ and the proofs for it are analogous. We
henceforth simply write $\xi$ for $\xi^-$ and $I$ for $I^-$.}

To demonstrate the required asymptotics for (\ref{stat1}) and (\ref{stat2}),
and to take advantage of the ideas just discussed, it is natural to
parametrize points in $Q$ as follows. First, assume without
loss of generality that $Q$ has Lebesgue measure equal to~$1$.
Let
$M(G)$ be the medial axis of $G$, that is, the closure of the
set of all points in $G$ with more than one closest point on
$\Gamma$. In general, $M(G)$ is a nonregular $(d-1)$-dimensional
surface with null $d$-dimensional Lebesgue measure. We refer the
reader to \cite{choi1997} for all matters concerning the theory of
medial axes.

Let $\Gamma_0 \subset\Gamma$ be the subset of $\Gamma$ where there
is no uniquely defined tangent plane. Then, ${\calH}^{d-1}(\Gamma_0)
= 0$ by assumption. Consider the subset $G_1$ of $G \setminus M(G)$
consisting of points $x$ such that there is a $\gamma\in\Gamma
\setminus\Gamma_0$ with $x - \gamma$ orthogonal to the tangent
plane at $\gamma$ and $\gamma:=\gamma(x)$ is the
boundary point closest to $x$.
Let $G_0$ be the
complement of $G_1$ with respect to $G \setminus M(G).$

For all $t > 0$, consider the level sets $\Gamma(t):= \{ x \in G
\backslash M(G)\dvtx d(x, \Gamma) = t \}$, where $d(x, \Gamma)$
denotes the distance between $x$ and $\Gamma$. Except possibly on a
Lebesgue null subset of $G_1$, we may uniquely parameterize points
$x \in G_1$ as $x(\gamma, t)$, where $\gamma$ is the boundary point
closest to $x$ and $t := \Vert x - \gamma\Vert$. Set $\Gamma(0) = \Gamma$.
Let $\Gamma_1(t):= \Gamma(t) \cap G_1$ and $\Gamma_0(t):= \Gamma(t)
\cap G_0$.
For each $\gamma\in
\Gamma\setminus\Gamma_0$, let $R_\gamma$ be the distance between
$\gamma$ and $M(G)$, measured along the orthogonal to the tangent
plane to $\Gamma$ at $\gamma$. Define $D := \sup_{\gamma\in\Gamma}
R_\gamma. $

For any integrable function $g\dvtx Q\to\R$, an application of the
co-area formula (see Theorem 3.2.12 and Lemma 3.2.34 in
\cite{federer1969}) for the distance function $f(x):= d(x,
\Gamma)$ gives
\[
\int_G g(x)\,dx = \int_0^D\!\!\! \int_{\Gamma(t)} {g(x)}{\calH}^{d-1}(dx)\,dt
\vadjust{\goodbreak}\]
since $|\nabla(f(x))| \equiv1$ a.e. \vadjust{\goodbreak}This yields the \textit{scaled
volume identity}

%e30 ###
%
\begin{eqnarray}\label{changeofvar1}
n^{1/d}\int_G g(x)\,dx &=& \int_0^{Dn^{1/d}}
\biggl[ \int_{ \Gamma_1(tn^{-1/d})} {g(x(\gamma,tn^{-1/d})) }
{\calH}^{d-1}(dx)\biggr]\,dt\nonumber\\ [-8pt]\\ [-8pt]
&&{}+ \int_0^{Dn^{1/d}} \biggl[ \int_{ \Gamma_0(tn^{-1/d})} { g(x) }
{\calH}^{d-1}(dx)\biggr ]\,dt.\nonumber
\end{eqnarray}
We may similarly write
%e31 ###
%
\begin{eqnarray}\label{changeofvar2}
 n^{1/d}\int_{G^c} g(x)\,dx &=&
\int_0^{D'n^{1/d}} \biggl[ \int_{ \Gamma'_1(tn^{-1/d})}
g(x(\gamma,tn^{-1/d})) {\calH}^{d-1}(dx) \biggr]\,dt\nonumber\\ [-8pt]\\ [-8pt]
&&{}+ \int_0^{D'n^{1/d}} \biggl[ \int_{ \Gamma'_0(tn^{-1/d})} g(x)
{\calH}^{d-1}(dx)\biggr]\,dt,\nonumber
\end{eqnarray}
where $D', \Gamma'_1(\cdot)$ and $\Gamma'_0(\cdot)$ are now the
analogs of $D', \Gamma_1(\cdot)$ and $\Gamma_0(\cdot)$, respectively.

To simplify the notation and to set the stage for developing the
analog of (\ref{JY1}), we define, for all $x \in G_1$, the shifted
and $n^{1/d}$-dilated binomial point measures
%e32 ###
%
\begin{equation} \label{defP}
\calP_n(\gamma,t) := n^{1/d}\bigl( \calX_{n-1} - x(\gamma,tn^{-1/d})\bigr)\cup\{\0 \},
\end{equation}
together with the shifted and dilated bodies
%e33 ###
%
\begin{equation} \label{defG}
G_n(\gamma,t) := n^{1/d}\bigl( G - x(\gamma,tn^{-1/d})\bigr).
\end{equation}
We similarly define the shifted and $\la^{1/d}$-dilated Poisson
point measures
%e34 ###
%
\begin{equation} \label{defP1}
\calP_\la(\gamma,t) := \la^{1/d}\bigl( \calP_{\la} -
x(\gamma,t\la^{-1/d})\bigr )\cup\{\0 \}
\end{equation}
and the shifted and dilated bodies
%e35 ###
%
\begin{equation} \label{defG1}
G_\la(\gamma,t) := \la^{1/d}\bigl( G - x(\gamma,t\la^{-1/d})\bigr).
\end{equation}
More generally, for all $x \in G,$ we write
\[
\calP_n(x) := n^{1/d}( \calX_{n-1} - x)\cup\{\0 \}
\]
and similarly for $G_n(x), \calP_\la(x)$ and $G_\la(x)$. %[need this
%more general def. JY]

Roughly speaking, the dilated bodies $G_n(\gamma,t)$ are converging
locally around $x(\gamma,tn^{-1/d})$ to a half-space, whereas the
restrictions of the point processes $\calP_n(\gamma,t)$ to
$G_n(\gamma,t)$ are converging to a homogeneous point process on
this half-space [see (\ref{conti})].

The next result, a consequence of this observation, is similar to
(\ref{JY1}) and is the key to all that follows. It shows that for
all $\gamma\in\Gamma\setminus\Gamma_0$ and $t \in(0, \infty)$,
the Hausdorff measure of faces incident to $\0$ and belonging to the
inner sewing $\mathcal{S}^-(\calP_n(\gamma,t),\partial G_n(\gamma,t))$ converges in
distribution to the Hausdorff measure of faces incident to $\0$ and
belonging to the inner sewing $\mathcal{S}^-(\calP^{\0},\partial\HH^d_{-t})$.

\begin{lemma}\label{lemma1}
 For all $\gamma\in\Gamma\setminus\Gamma_0$ and
$t \in(0,
\infty),$ we have %as $n \to\infty$
%e36 ###
%
\begin{equation} \label{lim}
\hspace*{27pt}\xi(\0 , \calP_n(\gamma,t), G_n(\gamma,t)) \tod\xi(t); \qquad\xi(\0 ,
\calP_\la(\gamma,t), G_\la(\gamma,t)) \tod\xi(t).\hspace*{-12pt}
\end{equation}
\end{lemma}
\begin{pf} We will only prove the first result since the second
follows by identical methods. Using a rotation if necessary, we
will, without loss of generality, assume that the vector $\0 -
(\gamma, 0)$ is orthogonal to the boundary of the half-space
$\HH^d_{-t}$.

For an arbitrary point set $\Y\subset\R^d$ and body $D$, let
$\xi_B(y, \Y\cap D)$ be the sum of the Hausdorff measures of the
faces of the Delaunay triangulation of $\Y\cap D$ lying wholly
inside $D$ and incident to $y \in\Y$. Recall that, by construction,
the faces giving nonzero contribution to $\xi(\0 ,
\calP_n(\gamma,t), G_n(\gamma,t))$ have vertices belonging to
$G_n(\gamma,t)$. Since the boundary of $G_n(\gamma,t)$ is
differentiable, it follows for large enough $n$ that these faces
will eventually belong entirely to $G_n(\gamma,t)$. In other words,
%e37 ###
%
\begin{equation} \label{contin}
\hspace*{36pt}\lim_{n \to\infty} \bigl|\xi(\0 , \calP_n(\gamma,t), G_n(\gamma,t)) -
\xi_B\bigl(\0 , \calP_n(\gamma,t)\cap G_n(\gamma,t)\bigr)\bigr|= 0 \qquad\mbox{a.s.}\nonumber
\end{equation}

Since $\0$ is at a distance $t$ from the boundary of $G_n(\gamma,t),$ it
follows that the value of $\xi_B$ at $\0$ with respect to $\calH
\cap(\RR^{d-1} \cap(-\infty,t])$ is determined by the points of
$\calH\cap(\RR^{d-1} \cap(-\infty,t])$ in a ball centered at $\0$
and of radius $R:= \max(t, R_0),$ where $R_0$ is the stabilization radius
at $\0$ for the graph of the standard Delaunay triangulation of
$\calH\cup\0$. Thus, reasoning exactly as in
(\ref{metric}) and (\ref{JY1}), we have that $\calH\cap(\RR^{d-1}
\cap(-\infty,t])$ is a continuity point for the function
$g(\calA):= \xi_B(\0,\calA)$ with respect to the topology on $
\mathcal L$ induced by the metric $D$ at (\ref{metric}). Since
%e38 ###
%
\begin{equation} \label{conti}
\calP_n(\gamma,t) \cap
G_n(\gamma,t) \tod\calH\cap\bigl(\R^{d-1} \cap(-\infty, t]\bigr),
\end{equation}
it therefore follows by the continuous mapping theorem that
\[
\xi_B\bigl(\0, \calP_n(\gamma,t)\cap G_n(\gamma,t)\bigr)\tod\xi_B\bigl(\0,
\calH\cap\bigl(\RR^{d-1} \cap(-\infty,t]\bigr)\bigr) = \xi(t).
\]
Combining this with
(\ref{contin}), we have $\xi(\0 , \calP_n(\gamma,t), G_n(\gamma,t))
\tod\xi(t)$, which completes the proof.
\end{pf}

Given (\ref{lim}), one has convergence of the corresponding
expectations in (\ref{lim}), provided the random variables satisfy
the customary uniform integrability condition. This is the content
of the next lemma.

\begin{lemma}\label{lemma2}
 For all $\gamma\in\Gamma\setminus\Gamma_0$ and
all $t
> 0,$ we have
\[
\lim_{n\to\infty} \EE[\xi(\0 , \calP_n(\gamma,t), G_n(\gamma
,t))] =
\lim_{\la\to\infty} \EE[\xi(\0 , \calP_\la(\gamma,t),
G_\la(\gamma,t))]= \EE[\xi(t)].
\]
\end{lemma}

\begin{pf}
It is straightforward that the empty sphere criterion
characterizing Delaunay triangulations (see, e.g., Chapter 4 of
\cite{small1996}) implies that, uniformly in $t$, the volume of a
Delaunay simplex in $\calP_n(\gamma,t)$ has exponentially decaying
tails, this being equivalent to the tail probability that a
circumcircle (circumsphere in dimension greater than $2$) contains
no points from the binomial point process $ \calP_n(\gamma,t)$. It
follows that for
$ p > 0,$ there is a constant $C(p)$ such that
%e39 ###
%
\begin{equation} \label{unifint}
\sup_n \sup_{(\gamma,t)} \EE[|\xi(\0 , \calP_n(\gamma,t),
G_n(\gamma,t))|^p] \leq C(p),
\end{equation}
showing that
\[
\{\xi(\0 , \calP_n(\gamma,t), G_n(\gamma,t)), n \geq1\}
\]
are uniformly integrable. By Lemma \ref{lemma1}, we obtain the
desired convergence of $\EE[\xi(\0 , \calP_n(\gamma,t),
G_n(\gamma,t))]$. The convergence of $\EE[\xi(\0 ,
\calP_\la(\gamma,t), G_\la(\gamma,t))]$ follows using identical
methods.
\end{pf}

The next two lemmas are also consequences of exponential decay of
the volume of the circumspheres not containing points from $n^{1/d}
\calX_n$. The next result shows that the expectations appearing in
Lemma \ref{lemma2} are uniformly bounded, a technical result
foreshadowing the upcoming use of the dominated convergence theorem.

\begin{lemma}\label{lemma4}
 There is an integrable function $F\dvtx [0,\infty) \to
[0,\infty)$, with exponentially decaying tails, such that
\[
\sup_n \sup_{x \in\Gamma(tn^{-1/d})} \EE[\xi(\0 , \calP_n(x),
G_n(x))] \leq F(t)
\]
and
\[
\sup_\la\sup_{x \in\Gamma(t\la^{-1/d})} \EE[\xi(\0 ,
\calP_\la(x), G_\la(x))] \leq F(t).
\]
\end{lemma}

\begin{pf}
For all $x \in\Gamma(tn^{-1/d}),$ let
\[
E(x):=\{\0 \
\mbox{belongs to a face} f \mbox{ of a simplex in }
{\calD}(\calP_n(x)),f \cap\partial\HH_t^d \neq\varnothing\}.
\]
By definition, $\xi$ vanishes on $E^c(x)$. It is easy to see that
$P[E(x)]$ has exponentially decaying
tails in $t$, uniformly in $n$ and $\gamma$. %By (\ref{unifint}) we
The first result follows by considering $\xi(\0 , \calP_n(x),
G_n(x)) {\bf1}(E(x)),$ and applying the Cauchy--Schwarz inequality
and the bound (\ref{unifint}) with $p = 2$. The second result
follows similarly.
\end{pf}

Our last lemma provides a high probability bound on the diameter of
Delaunay simplices, a result which follows immediately from the
exponential decay of the volume of spheres not containing points
from $\calX_n$ or $\calP_\lambda$. In other words, making use of
the bounds $P[\calP_\lambda\cap B_r(x) = \varnothing] = \exp(-\la r^d
v_d)$, where $v_d$ is the volume of the unit radius $d$-dimensional
ball, as well as the bounds $P[\X_n \cap B_r(x) = \varnothing] = (1 -
r^d v_d)^n$, we obtain the following result.

\begin{lemma}\label{lemma3}
For any constant $A > 0$, there is a constant $C > 0$
such that, with probability exceeding $1 - n^{-A}$, the diameter of
all Delaunay simplices in ${\mathcal D}(\X_n)$ is bounded by $C (\log
n/n)^{1/d}$. Thus, with probability exceeding $1 - n^{-A}$, we
have
\[
\sup_{i \leq n} \xi(X_i,\calX_n, G) \leq C (\log n/n)^{(d-1)/d}
\]
and
\[
\sup_{x \in\calP_\lambda} \xi(x,\calP_\lambda, G) \leq C
(\log\la/\la)^{(d-1)/d}.
\]
\end{lemma}

We now have all of the ingredients needed to prove Theorem
\ref{expectation}.

\begin{pf*}{Proof of Theorem \ref{expectation}}
The proof has two parts: the first
shows expectation convergence and the second uses concentration
inequalities to deduce the a.s. convergence.

\textit{Part} I. We use the above lemmas (binomial input) to first show
that\break
$\lim_{n\to\infty} \EE[I_{n}(h,\Gamma) ] = \alpha_d\int_\Gamma h\,d\Gamma$ for $h$ inner continuous, recalling that we
notationally simplify $\xi^-$ to
$\xi$ and $I_n^-$ to $I_n$. We omit the proof of $\lim_{\la\to
\infty} \EE[{\mathcal I}_{\la}^-(h,\break\Gamma) ] = \alpha_d\int_\Gamma h\,d\Gamma$ as it follows verbatim, using instead the
Poisson input versions
of the above lemmas.

Note that $\xi(X_1,\calX_n, G) = 0$ if $X_1\notin G$. Using
(\ref{scaling}) and conditioning on $X_1$, we have
\begin{eqnarray*}
\EE[I_{n}(h,\Gamma)] &=& n \EE[h(X_1) \xi(X_1, \calX_n, G)]\\
&=&n^{1/d} \EE\bigl[h(X_1) \xi\bigl(\0, n^{1/d}(\calX_n-X_1),
n^{1/d}(G-X_1)\bigr)\bigr] \\
&=&n^{1/d}\int_G h(x) \EE\bigl[\xi\bigl(\0, n^{1/d}(\calX_{n-1}-x)\cup
\{\0\}, n^{1/d}(G-x)\bigr)\bigr]\,dx.
\end{eqnarray*}
Using the scaled volume identity (\ref{changeofvar1}) and putting,
for all $x \in G$,
\[
G_n(x) := h(x) \EE\bigl[\xi\bigl(\0, n^{1/d}(\X_{n-1} - x) \cup\{\0 \},
n^{1/d}(G - x)\bigr)\bigr],
\]
we get that
%e40 ###
%
\begin{eqnarray} \label{2int}
\EE[I_{n}(h,\Gamma)] &=&\int_0^{Dn^{1/d}}
\int_{\Gamma_1(tn^{-1/d}) } G_n(x) {\calH}^{d-1}(dx)\,dt\nonumber\\
[-8pt]\\ [-8pt]
&&{}+ \int_0^{Dn^{1/d}} \int_{\Gamma_0(tn^{-1/d}) } G_n(x)
{\calH}^{d-1}(dx)\,dt.\nonumber
\end{eqnarray}
The proof of expectation convergence will be complete once we show
that the two integrals in (\ref{2int}) converge to $\alpha_d \
\int_\Gamma h(\gamma) {\calH}^{d-1}(d\gamma)$ and zero,
respectively.

We first consider the first integral in (\ref{2int}). For fixed $t
> 0$ and all $n,$ there is an a.e. $C^1$ mapping $f_n\dvtx \Gamma\to
\Gamma_1(tn^{-1/d})$
of $\Gamma$ onto the level set $\Gamma_1(tn^{-1/d})$.
To prepare for an application of the dominated
convergence theorem, we next show, for each $t
> 0,$ that as $n \to\infty,$
%e41 ###
%
\begin{equation} \label{limit1}
\int_{\Gamma_1(tn^{-1/d})} G_n(x) {\calH}^{d-1} (dx) \to
\int_{\Gamma} h(\gamma) \EE[(\xi(t))]{\calH}^{d-1}
(d\gamma).
\end{equation}
To see this, write the difference of the
integrals in (\ref{limit1}) as the sum of
%e42 ###
%
\begin{equation} \label{limit1a}
\int_{\Gamma_1(tn^{-1/d})} G_n(x) {\calH}^{d-1} (dx) - \int_{\Gamma}
G_n(f_n(\gamma)) {\calH}^{d-1}(d\gamma)
\end{equation}
and
%e43 ###
%
\begin{equation}\label{limit1b}
\int_{\Gamma} G_n(f_n(\gamma)) {\calH}^{d-1}
(d\gamma) - \int_{\Gamma} h(\gamma) \EE[(\xi(t))] {\calH}^{d-1}
(d\gamma).
\end{equation}
By a change of variables, the difference
(\ref{limit1a}) equals
\[
\int_\Gamma G_n(f_n(\gamma)) [ {\calH}^{d-1} (df_n(\gamma)) -
{\calH}^{d-1} (d\gamma)].
\]
By the a.e. smoothness of $\Gamma$,
we have $ {\calH}^{d-1} (df_n(\gamma)) = (1 + \epsilon_n(\gamma))
{\calH}^{d-1} (d\gamma)$ a.e., where $\epsilon_n(\gamma)$ goes to
zero as $n \to\infty$.
By (\ref{unifint}), we have
\[
\sup_n \sup_{\gamma\in\Gamma}|G_n(f_n(\gamma))| \leq C(1)
\Vert h\Vert_{\infty}
\]
and so the difference (\ref{limit1a}) goes to
zero. Next, consider the difference (\ref{limit1b}). Recalling that $t$ is fixed, for all $\gamma
\in\Gamma$ we write $f_n(\gamma) = x(\gamma, tn^{-1/d}).$ As $n \to
\infty,$ we have $ h(x(\gamma,tn^{-1/d})) \to h(x(\gamma,0)) =
h(\gamma),$ by inner continuity of $h$. Combining this with Lemma
\ref{lemma2}, we get, for all $x \in\Gamma_1(tn^{-1/d})$ with $x =
x(\gamma, tn^{-1/d})$, that $G_n(x)$ converges to $h(\gamma)
\EE[\xi(t)]$
and so, by the bounded convergence theorem, the difference
(\ref{limit1b}) goes to zero as $n \to\infty$. Thus, (\ref{limit1})
holds.

By Lemma \ref{lemma4}, we have, for fixed $t > 0,$
%e44 ###
%
\begin{equation} \label{Bound}
\hspace*{20pt}\int_{\Gamma_1(tn^{-1/d})} G_n(x) {\calH}^{d-1}(dx)
\leq\Vert h\Vert_{\infty} \times\sup_n (
{\calH}^{d-1}(\Gamma(tn^{-1/d}))) \times F(t).
\end{equation}
By the
boundedness of $h$ and ${\calH}^{d-1}(\Gamma)$, the right-hand side
of (\ref{Bound}) is integrable in $t$ and thus the dominated
convergence theorem implies
that
\[
\int_0^{Dn^{1/d}}\!\!\! \int_{\Gamma(tn^{-1/d}) } [
h(x(\gamma,tn^{-1/d})) \EE[\xi(\0,\calP_n(\gamma,t) ,
G_n(\gamma,t))] ] {\calH}^{d-1}(dx)\,dt
\]
converges to
\[
\int_0^{\infty}\!\!\!\int_{\Gamma} h(\gamma)
\EE[\xi(t)]
{\calH}^{d-1}(dx)\,dt = \int_0^{\infty}
\EE[\xi(t)]\,dt
\int_{\Gamma} h(\gamma) {\calH}^{d-1}(d\gamma),
\]
as desired.

Finally, we show that the second integral in (\ref{2int}) goes to
zero as $n \to\infty$. For fixed $t > 0,$ the inside integrals in
(\ref{2int}) are also
bounded by the right-hand side of (\ref{Bound}).
Since ${\calH}^{d-1}
(\Gamma_0(tn^{-1/d})) \to0$ as $n \to\infty$, the dominated
convergence theorem gives that the second integral in (\ref{2int})
converges to zero. This completes the proof of expectation
convergence.

\textit{Part} II. Next, we show a.s. convergence of $I_{n}(h,
\Gamma)$ and $I_{\la}(h, \Gamma)$. We will do this by appealing to a
variant of the Azuma--Hoeffding concentration inequality.
We only prove the convergence of $I_{n}(h, \Gamma)$ since the proof
of the convergence of $I_{\la}(h, \Gamma)$ is identical. Let $C_1$ be
a positive constant. For all $n,$
define the ``thickened boundary''
\[
G(n):= \{x \in G\dvtx d(x, \Gamma) \leq C_1(\log n/n)^{1/d}\}.
\]
By smoothness of $\Gamma,$ it follows that $v(n):=
\operatorname{volume}(G(n)) = O((\log n/n)^{1/d})$.
Define
%e45 ###
%
\begin{equation} \label{defI}
\hat{I'}_{n}(h,\Gamma):= \sum_{k=1}^{n v(n)} h(X_k) \xi\bigl(X_k,
{\calX}_{nv(n)}, G\bigr),
\end{equation}
where $X_k$, $k \geq1,$ belong to $G(n)$, ${\calX}_{nv(n)} :=
\{X_1,\ldots,X_{n v(n)}\}.$ In contrast to ${I}_{n}(h,\Gamma)$, note
that $\hat{I'}_{n}(h,\Gamma)$ contains a deterministic number of
nonzero terms and is therefore more amenable to analysis. Our goal
is to show that $\hat{I'}_{n}(h,\Gamma)$ well approximates
${I}_{n}(h,\Gamma)$ both a.s. and in $L^1$, and then to obtain
concentration results for $\hat{I'}_{n}(h,\Gamma)$. The proof of
a.s. convergence proceeds in the following four steps.

\textit{Step} (a). With high probability, the summands
$h(X_k) \xi(X_k, {\calX}_n, G)$ contributing a nonzero contribution
to ${I}_{n}(h,\Gamma)$ arise when $X_k$ belongs to the thickened
boundary $G(n)$. There are roughly $nv(n)$ such summands and thus
the convergence of ${I}_{n}(h,\Gamma)$ may be obtained by
restricting attention to the statistic $\hat{I'}_{n}(h,\Gamma)$
defined at (\ref{defI}). Our first goal is to make this precise.

The number of sample points in $\X_n$ belonging to $G(n)$ is a
binomial random variable $B(n, v(n))$. Relabeling, we may, without
loss of generality, assume that $X_1,\ldots,X_{B(n, v(n))}$ belong to
$G(n)$, where we suppress the dependency of $X_k$ on~$n$.

Define
\[
\hat{I}_{n}(h,\Gamma):= \sum_{k=1}^{B(n, v(n))} h(X_k)
\xi\bigl(X_k,{\calX}_{B(n, v(n))}, G\bigr),
\]
where ${\calX}_{B(n, v(n))} := \{X_1,\ldots,X_{B(n, v(n))}\}.$

By Lemma \ref{lemma3} and recalling the definition of $G(n)$, if
$C_1$ is large enough, then, with high probability, the simplices
defining the inner sewing of $G$ belong to $G(n)$ and thus it
follows that for any constant $A,$ there is a $C_1$ large enough such
that
%e46 ###
%
\begin{equation} \label{diff1}
P[\hat{I}_{n}(h,\Gamma) \neq{I}_{n}(h,\Gamma)] \leq n^{-A}.
\end{equation}
We will return to this bound in the sequel.

\textit{Step} (b). We next approximate
$\hat{I}_{n}(h,\Gamma)$ by $\hat{I'}_{n}(h,\Gamma)$. Given
$\hat{I}_{n}(h,\Gamma)$, we replace $B(n, v(n))$ by its mean,
which we assume, without loss of generality, to be integral (otherwise,
we use the integer part thereof). Observe that
\[
|\hat{I'}_{n}(h,\Gamma) -
\hat{I}_{n}(h,\Gamma)|
\]
is bounded by the product of four
factors:
\begin{enumerate}
\item[(i)] the difference between the cardinalities of the defining index
sets, namely $|B(n, v(n)) - n v(n)|$;
\item[(ii)] the maximal number $N$ of summands affected by either
deleting or inserting a single point into either ${\calX}_{B(n,
v(n))}$ or ${\calX}_{nv(n)}$;
\item[(iii)] $\sup_{k \leq nv(n)} \xi(X_k,{\calX}_{nv(n)}, G) +
\sup_{k
\leq B(n, v(n))} \xi(X_k,{\calX}_{B(n, v(n))}, G);$
\item[(iv)] the sup norm of $h$ on $G$, namely $\Vert h\Vert_{\infty}$.
\end{enumerate}
However, with high probability, we have these bounds:
\[
|B(n, v(n)) - n v(n)| \leq C_2 \log(n v(n)) (n v(n))^{1/2},
\]
$N \leq C_3 \log n$ (by Lemma \ref{lemma3}) and, by the analog of
Lemma \ref{lemma3},
\[
\sup_{k \leq nv(n)} \xi\bigl(X_k,{\calX}_{nv(n)}, G\bigr) + \sup_{k
\leq B(n, v(n))} \xi\bigl(X_k,{\calX}_{B(n, v(n)\bigr)}, G) \leq C_4( \log
n/n)^{(d-1)/d}.
\]
Here, and elsewhere, $C_1, C_2,\ldots$ denote generic
constants. We consequently obtain the high probability bound
\begin{eqnarray}\label{diff1a}
&&|\hat{I'}_{n}(h,\Gamma) - \hat{I}_{n}(h,\Gamma)|\nonumber\\
&&\qquad\leq C_2 [\log
(n v(n)) (n v(n))^{1/2}] [ C_3 \log n] \bigl[C_4(\log n/n)^{(d-1)/d}\bigr]
\Vert h\Vert_{\infty}\\
&&\qquad\leq C_5 (\log n)^3n^{-(d-1)/2d}.\nonumber
\end{eqnarray}

\textit{Step} (c). We combine (\ref{diff1}) and
(\ref{diff1a}), and take $A$ large enough in
(\ref{diff1}) to get
the high probability bound
\[
|\hat{I'}_{n}(h,\Gamma) - {I}_{n}(h,\Gamma)| \leq C_6 (\log n)^3
n^{-(d-1)/2d}.
\]
Since $ |\hat{I'}_{n}(h,\Gamma) -
{I}_{n}(h,\Gamma)|$ is deterministically bounded by a multiple of
$n$, it follows that
\[
\EE[|\hat{I'}_{n}(h,\Gamma) - {I}_{n}(h,\Gamma)|] \to0,
\]
whence
\[
\EE[\hat{I'}_{n}(h,\Gamma)] \to\alpha_d\int_\Gamma h(x)\,d\Gamma.
\]

It thus suffices to show that
%e48 ###
%
\begin{equation} \label{diff2}
|\hat{I'}_{n}(h,\Gamma) - \EE[\hat{I'}_{n}(h,\Gamma)]| \to0\qquad \mbox{a.s.}
\end{equation}

\textit{Step} (d). We complete the proof by showing
(\ref{diff2}). We do this by using a variant of the
Azuma--Hoeffding inequality, due to Chalker et al. \cite
{chalker1999}. Write $I(X_1,\ldots,X_{n v(n)})$ instead of
$\hat{I'}_{n}(h,\Gamma)$. Consider the martingale difference
representation
\[
I\bigl(X_1,\ldots,X_{nv(n)}\bigr) - \EE\bigl[I\bigl(X_1,\ldots,X_{nv(n)}\bigr)\bigr] =
\sum_{i=1}^{nv(n)}d_i,
\]
where $d_i:= \EE[ I(X_1,\ldots,X_{nv(n)}) |
\calF_i] - \EE[ I(X_1,\ldots,X_{nv(n)}) | \calF_{i-1}],$ here
$\calF_i$ being the $\sigma$-field generated by $X_1,\ldots,X_i$. Observe
that
\[
d_i:= \EE\bigl[ I\bigl(X_1,\ldots,X_{nv(n)}\bigr) | \calF_i\bigr] - \EE\bigl[
I\bigl(X_1,\ldots,X'_i,\ldots,X_{nv(n)}\bigr) | \calF_{i}\bigr],
\]
where $X_i'$ signals
an independent copy of $X_i$. By the conditional Jensen inequality
and Lemma \ref{lemma3}, it follows that
\begin{eqnarray*}
|d_i| &\leq&\EE\bigl[ \bigl|I\bigl(X_1,\ldots,X_{nv(n)}\bigr) -
I\bigl(X_1,\ldots,X'_i,\ldots,X_{nv(n)}\bigr)\bigr| \calF_{i}\bigr]\\
&\leq& C_7 (\log
n/n)^{(d-1)/d}
\end{eqnarray*}
holds on a high probability set, that is, for all $A >0,$ there is a
$C_7$ such that $P[|d_i| \geq C_7 (\log n/n)^{(d-1)/d}] \leq
n^{-A}$. If the $(d_i)_i$ were uniformly bounded in sup norm by
$o(1)$, then we could use the Azuma--Hoeffding inequality. Since this
is not the case, we use the following variant (see Lemma 1 of
\cite{chalker1999}), valid for all positive scalars $w_i, i \geq
1$:
%e49 ###
%
\begin{eqnarray}\label{Azuma}
P\Biggl[ \Biggl| \sum_{i=1}^{nv(n)} d_i \Biggr| > t \Biggr] &\leq&
2 \exp\biggl( \frac{-t^2}{32 \sum_{i = 1}^{nv(n)} w_i^2 }\biggr )
\nonumber\\ [-8pt]\\ [-8pt]
&&{} + \Bigl(1 + 2 t^{-1} \sup_{i \leq{nv(n)}} \Vert d_i \Vert_{\infty}\Bigr)
\sum_{i=1}^{nv(n)} P [ | d_i | > w_i ].\nonumber
\end{eqnarray}
Let $w_i = C_7 (\log n/n)^{(d-1)/d}$. We have $\sum_{i = 1}^{nv(n)}
w_i^2= C_8 (\log n)^{2-(1/d)}n^{-(d-1)/d},$ showing that the first term
in (\ref{Azuma}) is summable in $n$. Since $\Vert d_i \Vert_{\infty} \leq
C_9$ and $P[|d_i| \geq w_i] \leq n^{-A}$,
the second term in (\ref{Azuma}) is summable in $n$. Since $t$ is arbitrary,
this gives (\ref{diff2}), as desired.
\end{pf*}
\begin{pf*}{Proof of Theorem \ref{variance}}
We will follow ideas given in
\cite{schreiber2008}, which also involves functionals $\xi$ whose
expectations decay exponentially fast with the distance to the
boundary. For completeness, we provide the details, following in part
\cite{penrose2007EJP} and \cite{baryshnikov2009}.

Our goal is to show that
\[
\lim_{\la\to\infty} \la^{(d-1)/d} \Var[I_{\la}(h,\Gamma) ] =
V_d\int_\Gamma h^2(\gamma){\mathcal H}^{d-1}(d\gamma),
\]
where $I_{\la}(h,\Gamma)$ and $V_d$ are defined at (\ref{Pois}) and
(\ref{Vd}), respectively.

Let $\xi_\la(x,\calP_{\la}, G):= \xi(\la^{1/d}x,
\la^{1/d}\calP_{\la}, \la^{1/d}G).$ By scaling (\ref{scaling}), we
have
\[
\la^{(d-1)/d} \Var[I_{\la}(h,\Gamma)] = \la^{-(d-1)/d} \Var
\biggl[\sum_{x \in\calP_{\la} } h(x) \xi_\la(x,\calP_{\la},
G)\biggr].
\]
On the other hand, Campbell's theorem (see Chapter 13 of
\cite{daley2008}) gives
%e50 ###
%
\begin{eqnarray}\label{var1}
&&\la^{-(d-1)/d} \Var\biggl[\sum_{x \in\calP_{\la} } h(x) \xi_\la(x,\calP_{\la},
G)\biggr]\nonumber\\
&&\qquad= \la^{(d + 1)/d} \int_G \int_{\RR^d} [\cdots] h(x)h(y)\,dy\,dx\\
&&\qquad\quad{} + \la^{1/d} \int_G \EE[\xi^2_\la(x,\calP_{\la}, G)] h^2(x)\,dx,\nonumber
\end{eqnarray}
where
\[
[\cdots] := \EE[\xi_\la(x,\calP_{\la} \cup y, G) \xi_\la(y,\calP_{\la}
\cup x, G)] - \EE[\xi_\la(x,\calP_{\la}, G)]
\EE[\xi_\la(y,\calP_{\la}, G)].
\]
As in \cite{penrose2007EJP}, in the double integral in (\ref{var1}),
we put $y = x + \la^{-1/d}z$, thus giving
%e51 ###
%
\begin{eqnarray}\label{var2}
 &&\la^{(d + 1)/d} \int_G \int_{\RR^d} [\cdots] h(x)h(y)\,dy\,dx\nonumber\\ [-8pt]\\
 [-8pt]
 &&\qquad=\la^{1/d} \int_G \int_{\RR^d} F_\la(z,x) h(x)h(x +
 \la^{-1/d}z)\,dz\,dx,\nonumber
\end{eqnarray}
where
\begin{eqnarray*}
F_\la(z,x) &:=& \EE[\xi_\la(x,\calP_{\la} \cup\{x + \la^{-1/d}z\}
, G) \xi_\la(x + \la^{-1/d}z,\calP_{\la}
\cup x, G)]\\
&&{}- \EE[\xi_\la(x,\calP_{\la}, G)] \EE[\xi_\la(x +
\la^{-1/d}z,\calP_{\la}, G)]\hspace*{15pt}
\end{eqnarray*}
and where we adopt the convention that $\xi(x, \Y, G)$ is short for
$\xi(x, \Y\cup x, G)$ when $x$ is not in $\Y$. By the definition of
$\xi_\la$ and by translation invariance, we obtain that $F_\la(z,x)$
is equal to
\begin{eqnarray*}
&&\EE\bigl[\xi\bigl(\0,\la^{1/d}(\calP_{\la} - x) \cup\{z\}, \la^{1/d}(G -
x)\bigr) \xi\bigl(z,\la^{1/d}(\calP_{\la} - x) \cup\{\0\}, \la^{1/d}(G -
x)\bigr)\bigr]\\
&&\qquad{}-\EE\bigl[\xi\bigl(\0,\la^{1/d}(\calP_{\la} - x), \la^{1/d}(G - x)\bigr)\bigr]
\EE\bigl[\xi\bigl(z,\la^{1/d}(\calP_{\la} - x), \la^{1/d}(G - x)\bigr)\bigr].
\end{eqnarray*}

Write $x:=x(\gamma, t \lambda^{-1/d})$ and recall the definitions of
$\calP_\la(\gamma,t)$ and $G_\la(\gamma,t)$ from (\ref{defP1}) and
(\ref{defG1}), respectively, so that the above becomes %, that is
\begin{eqnarray*}
F_\la(z,x)&=& \EE\bigl[\xi\bigl(\0,\calP_\la(\gamma,t) \cup\{z\}, \la
^{1/d}(G -x)\bigr) \xi\bigl(z,\calP_\la(\gamma,t), \la^{1/d}(G - x)\bigr)\bigr ]\\
&&{}- \EE\bigl[\xi\bigl(\0,\calP_\la(\gamma,t), \la^{1/d}(G - x)\bigr)\bigr]
\EE\bigl[\xi\bigl(z,\calP_\la(\gamma,t), \la^{1/d}(G - x)\bigr)\bigr].
\end{eqnarray*}

Recalling that ${\xi}_B(y, \Y\cap D)$ is the normalized Hausdorff
measure of the faces of the Delaunay triangulation of $\Y\cap D$
lying inside $D$ and incident to $y$, the above becomes
\begin{eqnarray*}
\hspace*{-0.5pt}F_\la(z,x)&=& \EE\bigl[\xi_B\bigl(\0,[\calP_\la(\gamma,t) \cup\{z\}] \cap
G_\la(\gamma,t)\bigr) \xi_B\bigl(z,[\calP_\la(\gamma,t) \cup\{z\}] \cap
G_\la(\gamma,t)\bigr)\bigr]\\
\hspace*{-0.5pt}&&{}-\EE\bigl[\xi_B\bigl(\0,\calP_\la(\gamma,t) \cap G_\la(\gamma,t)\bigr)\bigr]
\EE\bigl[\xi_B\bigl(z,[\calP_\la(\gamma,t) \cup\{z\}] \cap G_\la(\gamma,t)\bigr)\bigr].
\end{eqnarray*}

Next, we have a two-dimensional version of (\ref{lim}), namely
\begin{eqnarray*}
&&\bigl[\xi_B\bigl(\0,[\calP_\la(\gamma,t)\cup\{z\}] \cap G_\la(\gamma,t)\bigr),
\xi_B\bigl(z,[\calP_\la(\gamma,t) \cup\{z\}] \cap G_\la(\gamma,t)\bigr)\bigr]\\
&&\qquad\tod\bigl[\xi_B(\0, [\calP^{\0}\cup\{z\}] \cap\HH_t^d),\xi_B(z,
[\calP^{\0}\cup\{z\}]\cap\HH_t^d)\bigr],
\end{eqnarray*}
from which it follows from uniform integrability
that as $\la\to\infty,$
%e52 ###
%
\begin{equation} \label{var4}
F_\la(z,x) \to c_t(z),
\end{equation}
where $c_t(z)$ is as in (\ref{cdef}).

We now find the large $\la$ behavior of the integrals at
(\ref{var2}). Recalling \vadjust{\goodbreak}the scaled volume identity
(\ref{changeofvar1}) and recalling that $x=x(\gamma, t
\lambda^{-1/d})$, we get, after substitution, that
\begin{eqnarray} \label{var5}
&&\la^{(d + 1)/d} \int_G \int_{\RR^d} [\cdots] h(x)h(y)\,dy\,dx\nonumber \\
&&\qquad= \int_0^{D \la^{1/d}}\!\!\! \int_{\Gamma_1(t \la^{-1/d})}
\int_{\RR^d} J_\lambda(z,t,\gamma)\,dz\,{\calH}^{d-1}\,(dx)\,dt\\
&&\qquad\quad{}+ \int_0^{D\la^{1/d}}\!\!\! \int_{\Gamma_0(t \la^{-1/d})}
\int_{\RR^d} J_\lambda(z,t,\gamma)\,dz\,{\calH}^{d-1}(dx)\,dt,\nonumber
\end{eqnarray}
where
\[
J_\lambda(z,t,\gamma):= F_\la(z,x(\gamma, t \lambda^{-1/d}))
h(x(\gamma, t \lambda^{-1/d})) h\bigl(x(\gamma, t \lambda^{-1/d}) +
\la^{-1/d}z\bigr).
\]
Notice that $J_\lambda(z,t,\gamma)$ is dominated uniformly in $\la$
by a function $F(z,t,\gamma)$ decaying exponentially fast in $|z|$
and $t$. By the a.e. continuity of $h$ and the convergence
(\ref{var4}), the integrand in the first integral tends to
$c_t(z)h^2(\gamma)$ as $\la\to\infty$. Bounding the integrand by
$\Vert h\Vert_\infty^2 F(z,t,\gamma)$ and applying dominated convergence, we
obtain that as $\la\to\infty,$ the first integral in (\ref{var5})
converges to
\begin{equation} \label{eq1}
\int_0^{\infty} \!\!\!\int_{\Gamma} \int_{\RR^d} h^2(\gamma) c_t(z)\,dz\,{\calH}^{d-1}(d\gamma)\,dt.
\end{equation}
The second integral in (\ref{var5}) converges to zero, by same
methods used to show that the second integral in (\ref{2int}) goes
to zero. Indeed,
\[
\int_{\Gamma_0(t \la^{-1/d}) }
\int_{\RR^d} J_\lambda(z,t,\gamma)\,dz\,{\calH}^{d-1}(dx)
\]
is bounded by an integrable function of $t$ which is going to zero
as $\la\to\infty$ since $\int_{\RR^d} J_\lambda(z,t,\gamma)\,dz$ are bounded uniformly in $\gamma$ and $t,$ and
${\calH}^{d-1}(\Gamma_0(t \la^{-1/d})) \to0$.

On the other hand, the single integral at (\ref{var1}) satisfies the
identity
\[
\la^{1/d} \int_G \EE[\xi^2_\la(x,\calP_{\la}, G)] h^2(x)\,dx =
\la^{1/d} \int_G \EE[\xi^2_\la(x,\calP_{\la}, G)] h^2(x)\,dx,
\]
which, as $\la\to\infty$, tends to
%e55 ###
%
\begin{equation} \label{eq2}
\int_0^{\infty} \int_{\Gamma} \EE[\xi^2(t)]\,dt\,h^2(\gamma) {\calH}^{d-1}(d\gamma)\,dt.
\end{equation}
Combining (\ref{eq1}) and (\ref{eq2}), and recalling the definition
of $V_d$ at (\ref{Vd}), we obtain
\[
\lim_{\la\to\infty} \la^{(d-1)/d} \Var[I_{\la}(h,\Gamma) ] =
V_d \int_\Gamma h^2(\gamma){\calH}^{d-1}(d\gamma).
\]
This completes the proof of Theorem \ref{variance}.
\end{pf*}

\section*{Acknowledgments}
Part of this work was done while R. Jim\'enez was visiting the
Department of Mathematics of Lehigh University---he wishes to thank the
faculty and staff of this department for their hospitality. The authors
also wish to express their gratitude to Jan Rataj and Joe Fu, as well
as to the anonymous referees
for their comments which led to an improved exposition.
%

%suskaldyti doi

%

%

\printaddresses

\end{document}